\documentclass{amsart}
\usepackage{verbatim,amsmath,amssymb,amscd,color,graphics,axodraw}
\usepackage{myyoungtab}

\usepackage[all]{xy} 
\CompileMatrices

\newcommand{\aal}{\mathbf{\stackrel{\leftarrow}{a}}}
\newcommand{\aar}{\mathbf{\stackrel{\rightarrow}{a}}}

\newcommand{\bij}{\Phi}
\newcommand{\Bt}{\tilde{B}}

\newcommand{\gc}{\geh_0}
\newcommand{\geh}{\mathfrak{g}}

\newcommand{\inner}[2]{\langle #1\,,\,#2\rangle}
\newcommand{\is}{\mathrm{inner}}

\newcommand{\la}{\lambda}
\newcommand{\La}{\Lambda}

\newcommand{\om}{\omega}
\newcommand{\os}{\mathrm{outer}}

\newcommand{\sigD}{\mathfrak{S}}
\newcommand{\ve}{\varepsilon}

\newcommand{\wt}{\mathrm{wt}\,}
\newcommand{\mb}{\mbox{}}
\newcommand{\bluem}{\textcolor{blue}{-}}
\newcommand{\redm}{\textcolor{red}{-}}

\newcommand\cb{{\bf c}}
\newcommand\cd{\cdots}
\newcommand\et[1]{\tilde{e}_{#1}}
\newcommand\ft[1]{\tilde{f}_{#1}}
\newcommand\Kz{K_{\Z}}
\newcommand\ol{\overline}
\newcommand\ot{\otimes}
\newcommand\Q{{\mathbb Q}}
\newcommand\Uq{U_q(\geh)}
\newcommand\Uqp{U'_q(\geh)}
\newcommand\Z{{\mathbb Z}}

\newcommand{\hdom}{{\begin{picture}(8,4)\multiput(0,0)(4,0){3}{\line(0,1){4}}
\multiput(0,0)(0,4){2}{\line(1,0){8}}\end{picture}}}
\newcommand{\vdom}{{\begin{picture}(4,8)\multiput(0,0)(0,4){3}{\line(1,0){4}}
\multiput(0,0)(4,0){2}{\line(0,1){8}}\end{picture}}}
\newcommand{\cell}{{\begin{picture}(4,4)\multiput(0,0)(4,0){2}{\line(0,1){4}}
\multiput(0,0)(0,4){2}{\line(1,0){4}}\end{picture}}}

\numberwithin{equation}{section}

\newtheorem{theorem}{Theorem}
\newtheorem{proposition}[theorem]{Proposition}
\newtheorem{lemma}[theorem]{Lemma}

\theoremstyle{definition}

\newtheorem{definition}{Definition}
\newtheorem{remark}{Remark}

\numberwithin{theorem}{section}
\numberwithin{definition}{section}
\numberwithin{remark}{section}

\begin{document}

\title[Existence of KR crystals]
{Existence of Kirillov--Reshetikhin crystals for nonexceptional types}

\author[Masato Okado]{Masato Okado}
\address{Department of Mathematical Science,
Graduate School of Engineering Science, Osaka University,
Toyonaka, Osaka 560-8531, Japan}
\email{okado@sigmath.es.osaka-u.ac.jp}

\author[Anne Schilling]{Anne Schilling}
\address{Department of Mathematics, University of California, One Shields
Avenue, Davis, CA 95616-8633, U.S.A.}
\email{anne@math.ucdavis.edu}
\urladdr{http://www.math.ucdavis.edu/\~{}anne}
\thanks{\textit{Date:} October 2008}

\dedicatory{Dedicated to Professor Masaki Kashiwara on his sixtieth birthday}

\begin{abstract}
Using the methods of~\cite{KMN2} and recent results on the characters of Kirillov--Reshetikhin 
modules \cite{H1,H2,N2}, the existence of Kirillov--Reshetikhin crystals $B^{r,s}$
is established for all nonexceptional affine types. 
We also prove that the crystals $B^{r,s}$ of type $B_n^{(1)}$, $D_n^{(1)}$, and $A_{2n-1}^{(2)}$
are isomorphic to the combinatorial crystals of~\cite{S:2007} for $r$ not a spin node.
\end{abstract}

\maketitle

\section{Introduction}
The theory of crystal bases by Kashiwara~\cite{K1} provides a remarkably powerful tool to study the 
representations of quantum algebras $U_q(\geh)$. For instance, the calculation of tensor product multiplicities
reduces to counting the number of crystal elements having certain properties. Although crystal bases
are bases at $q=0$, one can ``melt" them to get actual bases, called global crystal bases, for integrable
highest weight representations of $\Uq$. It turns out that the global crystal basis agrees with Lusztig's
canonical basis~\cite{L:Book}, and it has many applications in representation theory.

The main focus of this paper are affine finite crystals, that is, crystal bases of finite-dimensional
modules for quantum groups 
corresponding to affine Kac--Moody algebras $\geh$. These crystal bases were first developed by
Kang et al.~\cite{KMN1,KMN2}, where it was also shown that integrable highest-weight
$U_q(\geh)$-modules of arbitrary level can be realized as semi-infinite tensor
products of perfect crystals. This is known as the path realization.
Many perfect crystals were proven to exist and explicitly constructed in~\cite{KMN2}.

Irreducible finite-dimensional $U_q'(\geh)$-modules were classified by Chari
and Pressley~\cite{CP,CP:1998} in terms of Drinfeld polynomials. It was conjectured 
by Hatayama et al.~\cite{HKOTT,HKOTY} that a certain subset of such modules known 
as Kirillov--Reshetikhin (KR) modules $W^{(r)}_s$ have a crystal basis $B^{r,s}$.
Here the index $r$ corresponds to a node of the Dynkin diagram of $\geh$ except the prescribed
$0$ and $s$ is an arbitrary positive integer.
This conjecture was confirmed in many instances~\cite{BFKL,KMN1,KMN2,K3,
Koga:1999,NS:2006,Y:1998}, but a proof for general $r$ and $s$ has not been available except 
type $A_n^{(1)}$ in \cite{KMN2}. Only recently the existence proof was completed 
in~\cite{O:2006} for type $D_n^{(1)}$.
Using the methods of~\cite{KMN2} and recent results on the characters of KR modules~\cite{H1,H2,N2}, 
we establish the existence of Kirillov--Reshetikhin crystals $B^{r,s}$
for all nonexceptional affine types in this paper:
\begin{theorem} \label{thm:main}
The Kirillov-Reshetikhin module $W^{(r)}_s$ associated to any nonexceptional affine
Kac--Moody algebra has a crystal basis $B^{r,s}$.
\end{theorem}

In addition we prove that for type $B_n^{(1)}$, $D_n^{(1)}$, and $A_{2n-1}^{(2)}$
these crystals coincide with the combinatorial crystals of~\cite{S:2007,St:2006}.
Throughout the paper we denote by $B^{r,s}$ the KR crystal associated with the KR module
$W_s^{(r)}$. The combinatorial crystal of~\cite{S:2007} is called $\Bt^{r,s}$.
Our second main result is the following theorem:
\begin{theorem}  \label{thm:equiv}
For $1\le r\le n-2$ for type $D_n^{(1)}$, $1\le r\le n-1$ for type $B_n^{(1)}$, $1\le r\le n$ for 
type $A_{2n-1}^{(2)}$ and $s\in\Z_{>0}$, the crystals $B^{r,s}$ and
$\Bt^{r,s}$ are isomorphic.
\end{theorem}

The key to the proof of Theorem~\ref{thm:main} is Proposition~\ref{prop:key} below, which is 
due to Kang et al.~\cite{KMN2} and states that a finite-dimensional $U_q'(\geh)$-module having 
a prepolarization and certain $\Z$-form has a crystal basis if the dimensions of some particular
weight spaces are not greater than the weight multiplicities of a fixed module and the values of
the prepolarization of certain vectors in the module have some special properties. Using the
fusion construction it is established that the KR modules have a prepolarization and $\Z$-form.
The requirements on the dimensions follow from recent results by Nakajima~\cite{N2}
and Hernandez~\cite{H1,H2}. Necessary values of the prepolarization are calculated 
explicitly in Propositions~\ref{prop:norm DBA}, \ref{prop:norm C}, and \ref{prop:norm AD}.

The isomorphism between the KR crystal $B^{r,s}$ and the combinatorial crystal
$\Bt^{r,s}$ is established by showing that isomorphisms as crystals with index
sets $\{1,2,3,\ldots,n\}$ and $\{0,2,3,\ldots,n\}$ already uniquely determine the whole crystal.

Before presenting our results, let us offer some speculations on combinatorial realizations
for the KR crystals. For type $A_n^{(1)}$ the crystals $B^{r,s}$ were constructed combinatorially 
by Shimozono~\cite{Sh:2002} using the promotion operator. The promotion operator $\mathrm{pr}$
is the crystal analogue of the Dynkin diagram automorphism that maps node $i$ to node $i+1$ 
modulo $n+1$. The affine crystal operator $\ft{0}$ is then given by $\ft{0}=\mathrm{pr}^{-1}\circ \ft{1}
\circ \mathrm{pr}$. Similarly, the main tool used in~\cite{S:2007} to construct the combinatorial 
crystals $\Bt^{r,s}$ of type $B_n^{(1)}$, $D_n^{(1)}$, and $A_{2n-1}^{(2)}$ is the crystal analogue of the
Dynkin diagram automorphism that interchanges nodes 0 and 1.
For type $C_n^{(1)}$ and $D_{n+1}^{(2)}$, there exists a Dynkin diagram automorphism
$i\mapsto n-i$. It is our intention to exploit this symmetry to construct $\Bt^{r,s}$
of type $C_n^{(1)}$ and $D_{n+1}^{(2)}$ explicitly in a future publication.
For type $A_{2n}^{(2)}$ no Dynkin diagram automorphism exists. However, it should
still be possible to construct these crystals by looking at the $\{1,2,\ldots,n\}$ and 
$\{0,1,2,\ldots,n-1\}$ subcrystals as was done for $r=1$ in~\cite{KMN2}.
Realizations of $B^{r,s}$ as virtual crystals were given in~\cite{OSS:2003,OSS:2003a}.

The paper is organized as follows. In Section~\ref{sec:fr} we review necessary background
on the quantum algebra $U_q'(\geh)$ and the fundamental representations. In particular
we review Proposition~\ref{prop:key} of~\cite{KMN2} which provides a criterion for the
existence of a crystal pseudobase. In Section~\ref{sec:KR} we define KR modules
by the fusion construction and show that these modules have a prepolarization.
This reduces the existence proof for KR crystals to conditions stated in 
Proposition~\ref{prop:main}. These conditions are checked explicitly in 
Section~\ref{sec:existence} for the various types to prove Theorem~\ref{thm:main}.
In Section~\ref{sec:combinatorial KR} we review the combinatorial construction
of the crystals $\Bt^{r,s}$ of types $B_n^{(1)}$, $D_n^{(1)}$, and $A_{2n-1}^{(2)}$ and prove
in Section~\ref{sec:equivalence} that they are isomorphic to $B^{r,s}$, thereby establishing
Theorem~\ref{thm:equiv}.

\subsection*{Acknowledgments}

M.O. thanks Masaki Kashiwara for letting him know that the irreducibility 
of $W^{(r)}_s$ in Proposition~\ref{prop:KR irred} follows from his result in~\cite{K3},
and Hiraku Nakajima for attracting his attention to references on the polarization
of $V(\varpi_r)$. The authors would also like to thank David Hernandez for helpful
correspondences. M.O. is partially supported by Grant-in-Aid for Scientific Research (C) 18540030,
Japan Society for the Promotion of Science.
A.S. is partially supported by NSF grants DMS-0501101, DMS-0652641, and
DMS-0652652.

\subsection*{Note added after publication}
After publication we noticed some errors and omissions in our paper, which are
corrected in the ``Erratum'' in Appendix~\ref{erratum} at the end of the paper.
Also, Theorem~\ref{thm:equiv} has now been extended to all nonexceptional types
in~\cite{FOS:2008}.

\section{Quantum affine algebra $\Uqp$ and fundamental representations}
\label{sec:fr}

\subsection{Quantum affine algebra}

Let $\geh$ be an affine Kac-Moody algebra and $\Uq$ the quantum affine algebra associated 
to $\geh$. In this section $\geh$ can be any affine algebra. For the notation of $\geh$ or
$\Uq$ we follow \cite{K3}. For instance, $P$ is the weight lattice, $I$ is the index set of 
simple roots, and $\{\alpha_i\}_{i\in I}$ (resp. $\{h_i\}_{i\in I}$) is the set of simple
roots (resp. coroots). Let $(\;,\;)$ be the inner product on $P$ normalized by $(\delta,\la)
=\langle c,\la\rangle$ for any $\la\in P$ as in \cite{Kac}, where $c$ is the canonical central
element and $\delta$ is the generator of null roots.
We choose a positive integer $d$ such that $(\alpha_i,\alpha_i)/2\in
\Z d^{-1}$ for any $i\in I$ and set $q_s=q^{1/d}$. Then $\Uq$ is the associative algebra
over $\Q(q_s)$ with 1 generated by $e_i,f_i$ ($i\in I$), $q^h$ ($h\in d^{-1}P^*,P^*=
\mathrm{Hom}_{\Z}(P,\Z)$) with certain relations. By convention, we set $q_i=
q^{(\alpha_i,\alpha_i)/2},t_i=q_i^{h_i},[m]_i=(q_i^m-q_i^{-m})/(q_i-q_i^{-1}),
[n]_i!=\prod_{m=1}^n[m]_i,e_i^{(n)}=e_i^n/[n]_i!,f_i^{(n)}=f_i^n/[n]_i!$.

Let $\{\La_i\}_{i\in I}$ be the set of fundamental weights. Then we have 
$P=\bigoplus_i\Z\La_i\oplus\Z\delta$. We set 
\[
P_{cl}=P/\Z\delta.
\]
Similar to the quantum algebra $\Uq$ which is associated with $P$, we can also consider 
$\Uqp$, which is associated with $P_{cl}$, namely, the subalgebra of $\Uq$ generated by 
$e_i,f_i,q^h$ ($h\in d^{-1}(P_{cl})^\ast$).

Next we introduce two subalgebras (`$\Z$-forms') $\Uq_{\Kz}$ and $\Uq_{\Z}$ of $\Uq$.
Let $A$ be the subring of $\Q(q_s)$ consisting of 
rational functions without poles at $q_s=0$. We introduce the subalgebras $A_{\Z}$ and $\Kz$ 
of $\Q(q_s)$ by
\begin{align*}
A_{\Z}&=\{f(q_s)/g(q_s)\mid f(q_s),g(q_s)\in\Z[q_s],g(0)=1\},\\
\Kz&=A_{\Z}[q_s^{-1}].
\end{align*}
Then we have 
\[
\Kz\cap A=A_{\Z},\quad A_{\Z}/q_sA_{\Z}\simeq\Z.
\]
We then define $\Uq_{K_{\Z}}$ as the $K_{\Z}$-subalgebra of $\Uq$ generated by 
$e_i,f_i,q^h$ ($i\in I,h\in d^{-1}P^\ast$). $\Uq_{\Z}$ is defined as the $\Z[q_s,q_s^{-1}]$-subalgebra
of $\Uq$ generated by $e_i^{(n)},f_i^{(n)},{t_i\brace n}_i$ ($i\in I,n\in\Z_{>0}$) and $q^h$
($h\in d^{-1}P^\ast$). Here we have set ${x\brace n}_i=\prod_{k=1}^n(q_i^{1-k}x-q_i^{k-1}x^{-1})/[n]_i!$. 
$\Uq_{\Z}$ is a $\Z[q_s,q_s^{-1}]$-subalgebra of $\Uq_{\Kz}$. We can also introduce subalgebras 
$\Uqp_{\Kz}$ and $\Uqp_{\Z}$ by replacing $q^h$ ($h\in d^{-1}P^\ast$) with $q^h$ 
($h\in d^{-1}(P_{cl})^\ast$) in the generators.

We define a total order on $\Q(q_s)$ by
\[
f>g\text{ if and only if }f-g\in\bigsqcup_{n\in\Z}\{q_s^n(c+q_sA)\mid c>0\}
\]
and $f\ge g$ if $f>g$ or $f=g$. 

Let $M$ and $N$ be $\Uq$(or $\Uqp$)-modules. A bilinear form 
$(\;,\;):M\ot_{\Q(q_s)}N\rightarrow\Q(q_s)$ is called an admissible pairing if it satisfies
\begin{align}
(q^hu,v)&=(u,q^hv),\nonumber\\
(e_iu,v)&=(u,q_i^{-1}t_i^{-1}f_iv), \label{admissible}\\
(f_iu,v)&=(u,q_i^{-1}t_ie_iv),\nonumber
\end{align}
for all $u\in M$ and $v\in N$. Equation~\eqref{admissible} implies
\begin{equation} \label{admissible2}
(e_i^{(n)}u,v)=(u,q_i^{-n^2}t_i^{-n}f_i^{(n)}v),\quad
(f_i^{(n)}u,v)=(u,q_i^{-n^2}t_i^{n}e_i^{(n)}v).
\end{equation}
A symmetric bilinear form $(\;,\;)$ on $M$ is called a \textit{prepolarization} of $M$ if it
satisfies \eqref{admissible} for $u,v\in M$.
A prepolarization is called a \textit{polarization} if it is positive definite with respective to 
the order on $\Q(q_s)$.

\subsection{Criterion for the existence of a crystal pseudobase}

Here we recall the criterion for the existence of a crystal pseudobase given in \cite{KMN2}.
We do not review the notion of crystal bases, but refer the reader to \cite{K1}. We only note
that $q$ in the definition of crystal base in \cite{K1} should be replaced by $q_s$ according
to the normalization of the inner product $(\;,\;)$ on $P$. We say $(L,B)$ is a crystal pseudobase
of an integrable $\Uq$ (or $\Uqp$)-module $M$, if (i) $L$ is a crystal lattice of $M$, (ii) 
$B=B'\sqcup(-B')$ where $B'$ is a $\Q$-base of $L/q_sL$, (iii) $B=\bigsqcup_{\la\in P}B_{\la}$
where $B_{\la}=B\cap(L_{\la}/q_sL_{\la})$, (iv) $\et{i}B\subset B\sqcup\{0\}$,
$\ft{i}B\subset B\sqcup\{0\}$, and (v) for $b,b'\in B$, $b'=\ft{i}b$ if and only if $b=\et{i}b'$.
Note that only the condition (ii) is replaced from the definition of the crystal base.

Let $\geh_0$ be the finite-dimensional simple Lie algebra whose Dynkin 
diagram is obtained by removing the $0$-vertex from that of $\geh$. In this paper we specify the 
$0$-vertex as in \cite{Kac} and set $I_0=I\setminus\{0\}$. Let $\ol{P}_+$ be the set
of dominant integral weights of $\geh_0$ and $\ol{V}(\la)$ be the irreducible highest weight 
$U_q(\geh_0)$-module of highest weight $\la$ for $\la\in\ol{P}_+$. The following proposition is 
easily obtained by combining Proposition 2.6.1 and 2.6.2 of \cite{KMN2}.

\begin{proposition} \label{prop:key}
Let $M$ be a finite-dimensional integrable $\Uqp$-module. Let $(\;,\;)$ be a prepolarization
on $M$, and $M_{\Kz}$ a $\Uqp_{K_{\Z}}$-submodule of $M$ such that $(M_{\Kz},M_{\Kz})\subset\Kz$. 
Let $\la_1,\ldots,\la_m\in\ol{P}_+$, and assume that the following conditions hold:
\begin{equation} \label{char}
\dim M_{\la_k}\le\sum_{j=1}^m\dim \ol{V}(\la_j)_{\la_k}\mbox{ for }k=1,\ldots,m.
\end{equation}
\begin{align}
&\text{There exist $u_j\in(M_{\Kz})_{\la_j}$ $(j=1,\ldots,m)$ such that 
$(u_j,u_k)\in\delta_{jk}+q_sA$,} \label{norm}\\
&\text{and $(e_iu_j,e_iu_j)\in q_sq_i^{-2(1+\langle h_i,\la_j\rangle)}A$ for any $i\in I_0$.}
\nonumber
\end{align}
Set $L=\{u\in M\mid(u,u)\in A\}$ and set $B=\{b\in M_{\Kz}\cap L/M_{\Kz}\cap q_sL\mid(b,b)_0=1\}$.
Here $(\;,\;)_0$ is the $\Q$-valued symmetric bilinear form on $L/q_sL$ induced by $(\;,\;)$.
Then we have the following:
\begin{itemize}
\item[(i)] $(\;,\;)$ is a polarization on $M$.
\item[(ii)] $M\simeq\bigoplus_j \ol{V}(\la_j)$ as $U_q(\geh_0)$-modules.
\item[(iii)] $(L,B)$ is a crystal pseudobase of $M$.
\end{itemize}
\end{proposition}

\subsection{Fundamental representations}

For any $\la\in P$, Kashiwara defined a $\Uq$-module $V(\la)$ called extremal weight module \cite{K2}.
We briefly recall its definition. Let $W$ be the Weyl group associated to $\geh$ and $s_i$ the simple 
reflection for $\alpha_i$. Let $M$ be an integrable $\Uq$-module. A vector $u_{\la}$ of weight $\la\in P$ 
is called an extremal vector if there exists a set of vectors $\{u_{w\la}\}_{w\in W}$ satisfying 
\begin{align}
&u_{w\la}=u_{\la}\text{ for }w=e,\\
&\text{if }\langle h_i,w\la\rangle\ge0,\text{ then }e_iu_{w\la}=0\text{ and }
f_i^{(\langle h_i,w\la\rangle)}u_{w\la}=u_{s_iw\la},\\
&\text{if }\langle h_i,w\la\rangle\le0,\text{ then }f_iu_{w\la}=0\text{ and }
e_i^{(-\langle h_i,w\la\rangle)}u_{w\la}=u_{s_iw\la}.
\end{align}
Then $V(\la)$ is defined to be the $\Uq$-module generated by $u_{\la}$ with the defining relations that
$u_{\la}$ is an extremal vector. For our purpose, we only need $V(\la)$ when $\la=\varpi_r$ for 
$r\in I_0$, where $\varpi_r$ is a level $0$ fundamental weight
\begin{equation} \label{varpi}
\varpi_r=\La_r-\langle c,\La_r\rangle\La_0.
\end{equation}
Then the following facts are known.

\begin{proposition} \cite[Proposition 5.16]{K3}
\begin{itemize}
\item[(i)] $V(\varpi_r)$ is an irreducible integrable $\Uq$-module.
\item[(ii)] $\dim V(\varpi_r)_\mu<\infty$ for any $\mu\in P$.
\item[(iii)] $\dim V(\varpi_r)_\mu=1$ for any $\mu\in W\varpi_r$.
\item[(iv)] $\wt V(\varpi_r)$ is contained in the intersection of $\varpi_r+\sum_{i\in I}\Z\alpha_i$
	and the convex hull of $W\varpi_r$.
\item[(v)] $V(\varpi_r)$ has a global crystal base $(L(\varpi_r),B(\varpi_r))$.
\item[(vi)] Any integrable $\Uq$-module generated by an extremal weight vector of weight $\varpi_r$ is
	isomorphic to $V(\varpi_r)$.
\end{itemize}
\end{proposition}

Let $\la\in P^0=\{\la\in P \mid \inner{c}{\la}=0\}$. $V(\la)$ has a $\Uq_{\Z}$-submodule $V(\la)_{\Z}$. Let
$\{G(b)\}_{b\in B(\la)}$ stand for the global base of $V(\la)$. The following result was shown 
in~\cite{VV} for $\geh$ simply laced and $\la=\varpi_r$, in~\cite{N1} for $\geh$ simply laced and
$\la$ is arbitrary, and in~\cite{BN} for $\geh$ and $\la$ arbitrary.

\begin{proposition}\mbox{}
\begin{itemize}
\item[(i)] There exists a prepolarization $(\;,\;)$ on $V(\la)$.
\item[(ii)] $\{G(b)\}_{b\in B(\la)}$ is almost orthonormal with respect to $(\;,\;)$,
	that is, $(G(b),G(b'))\equiv\delta_{bb'}\text{ mod }q_s\Z[q_s]$.
\end{itemize}
\end{proposition}

Let $d_r$ be a positive integer such that
\[
\{k\in\Z \mid \varpi_r+k\delta\in W\varpi_r\}=\Z d_r.
\]
We note that $d_r=\max(1,(\alpha_r,\alpha_r)/2)$ except in the case $d_r=1$ when $\geh=A_{2n}^{(2)}$ 
and $r=n$. Then there exists a $\Uqp$-linear automorphism $z_r$ of $V(\varpi_r)$ of weight $d_r\delta$
sending $u_{\varpi_r}$ to $u_{\varpi_r+d_r\delta}$. Hence we can define a $\Uqp$-module $W(\varpi_r)$ by
\[
W(\varpi_r)=V(\varpi_r)/(z_r-1)V(\varpi_r).
\]
This module is called a fundamental representation. 

For a $\Uqp$-module $M$ let $M_{\mathrm{aff}}$ denote the $\Uqp$-module 
$\Q(q_s)[z,z^{-1}]\ot M$ with the actions of $e_i$ and $f_i$ by $z^{\delta_{i0}}\ot e_i$ and 
$z^{-\delta_{i0}}\ot f_i$. For $a\in\Q(q_s)$ we define the $\Uqp$-module $M_a$ by
$M_{\mathrm{aff}}/(z-a)M_{\mathrm{aff}}$.

\begin{proposition} \cite[Proposition 5.17]{K3}
\begin{itemize}
\item[(i)] $W(\varpi_r)$ is a finite-dimensional irreducible integrable $\Uqp$-module.
\item[(ii)] For any $\mu\in \wt V(\varpi_r)$, $W(\varpi_r)_{cl(\mu)}\simeq V(\varpi_r)_{\mu}$. Here 
	the map $cl$ stands for the canonical projection $P\longrightarrow P_{cl}$.
\item[(iii)] $\dim W(\varpi_r)_{cl(\mu)}=1$ for any $\mu\in W\varpi_r$.
\item[(iv)] $\wt W(\varpi_r)$ is contained in the intersection of $cl(\varpi_r+\sum_{i\in I}\Z\alpha_i)$
	and the convex hull of $W\,cl(\varpi_r)$.
\item[(v)] $W(\varpi_r)$ has a global crystal base.
\item[(vi)] Any irreducible finite-dimensional integrable $\Uqp$-module with $cl(\varpi_r)$ as an extremal
	weight is isomorphic to $W(\varpi_r)_a$ for some $a\in\Q(q_s)$.
\end{itemize}
\end{proposition}

We also need the following lemma that ensures the existence of the prepolarization on $W(\varpi_r)$.

\begin{lemma} \cite{VV,N1}
$(z_ru,z_rv)=(u,v)$ for $u,v\in V(\varpi_r)$.
\end{lemma}

\begin{remark} \label{rem:VV and N}
This lemma is given as Proposition 7.3 of \cite{VV} and also as Lemma 4.7 of \cite{N1}.
The lemmas or properties used to prove it hold for any affine algebra $\geh$.
\end{remark}

Summing up the above discussions we have

\begin{proposition}
The fundamental representation $W(\varpi_r)$ has the following properties:
\begin{itemize}
\item[(i)] $W(\varpi_r)$ has a polarization $(\;,\;)$.
\item[(ii)] There exists a $\Uqp_{\Z}$-submodule $W(\varpi_r)_{\Z}$ of $W(\varpi_r)$ such that 
\begin{equation*}
	(W(\varpi_r)_{\Z},W(\varpi_r)_{\Z})\subset\Z[q_s,q_s^{-1}].
\end{equation*}
\end{itemize}
\end{proposition}

Before finishing this section, let us mention the Drinfeld polynomials. It is known that irreducible
finite-dimensional $\Uqp$-modules are classified by $|I_0|$-tuple of polynomials $\{P_j(u)\}_{j\in I_0}$ 
whose constant terms are $1$. See e.g. \cite{CP}. The degree of $P_j$ is given by $\inner{\la}{h_j}$
where $\la$ is the highest weight of the corresponding module. Hence we have

\begin{lemma} \label{lem:Dpoly}
$W(\varpi_r)$ has the following Drinfeld polynomials
\[
P_r(u)=1-a^\dagger_ru,\quad P_j(u)=1\text{ for }j\neq r
\]
with some $a^\dagger_r\in\Q(q_s)$.
\end{lemma}

For types $A_n^{(1)},D_n^{(1)},E_{6,7,8}^{(1)}$ the explicit value of $a^\dagger_r$ is 
known~\cite[Remark 3.3]{N1}.

\section{KR modules and the existence of crystal bases}
\label{sec:KR}

\subsection{Fusion construction} \label{ss:fusion}

Let $V$ be a $\Uqp$-module. An $R$-matrix, denoted by $R(x,y)$, is an element of 
$\mbox{Hom}_{\Uqp[x^{\pm1},y^{\pm1}]}(V_x\ot V_y,V_y\ot V_x)$. For $V$ we assume the following:
\begin{align}
& \text{$V\ot V$ is irreducible}. \label{Virr}\\
& \text{There exists $\la_0\in P_{cl}$ such that 
 $\wt V\subset\la_0+\sum_{i\in I_0}\Z_{\le0}\alpha_i\text{ and }\dim V_{\la_0}=1$}. \label{la0}
\end{align}
Under these assumptions it is known (see e.g. \cite{KMN1}) that there
exists a unique $R$-matrix up to multiple of a scalar function of $x,y$.  
Take a nonzero vector $u_0$ from $V_{\la_0}$. We normalize $R(x,y)$ in such a way that
$R(x,y)(u_0\ot u_0)=u_0\ot u_0$. The normalized $R$-matrix is known to depend only on $x/y$.
Because of the normalization, some matrix elements of $R(x,y)$ may have zeros or poles as 
a function of $x/y$. At the points $x/y=x_0/y_0\in\Q(q_s)$ where there is no zero or pole, 
$R(x_0,y_0)$ is an isomorphism.

Next we review the fusion construction following section 3 of \cite{KMN2}.
Let $s$ be a positive integer and $\mathfrak{S}_s$ the $s$-th symmetric group.
Let $s_i$ be the simple reflection which interchanges $i$ and $i+1$, and 
let $\ell(w)$ be the length of $w\in\mathfrak{S}_s$. Let $R(x,y)$ denote the 
$R$-matrix for $V_x\ot V_y$. For any $w\in\mathfrak{S}_s$ we can construct a well-defined map 
$R_w(x_1,\ldots,x_s):V_{x_1}\ot\cd\ot V_{x_s}\rightarrow V_{x_{w(1)}}\ot\cd\ot V_{x_{w(s)}}$
by
\begin{align*}
R_1(x_1,\ldots,x_s)&=1,\\
R_{s_i}(x_1,\ldots,x_s)&=\left(\bigotimes_{j<i}\text{id}_{V_{x_j}}\right)\ot
R(x_i,x_{i+1})\ot\left(\bigotimes_{j>i+1}\text{id}_{V_{x_j}}\right),\\
R_{ww'}(x_1,\ldots,x_s)&=R_{w'}(x_{w(1)},\ldots,x_{w(s)})\circ R_w(x_1,\ldots,x_s)\\
&\hspace{1cm}\text{ for $w,w'$ such that $\ell(ww')=\ell(w)+\ell(w')$.}
\end{align*}
Fix $k\in d^{-1}\Z\setminus\{0\}$. Let us assume that 
\begin{equation} \label{pole}
\text{the normalized $R$-matrix $R(x,y)$ does not have a pole at $x/y=q^{2k}$.}
\end{equation}
For each $s\in\Z_{>0}$, we put
\begin{align*}
R_s=&R_{w_0}(q^{k(s-1)},q^{k(s-3)},\ldots,q^{-k(s-1)}):\\
&V_{q^{k(s-1)}}\ot V_{q^{k(s-3)}}\ot\cd\ot V_{q^{-k(s-1)}}\rightarrow
V_{q^{-k(s-1)}}\ot V_{q^{-k(s-3)}}\ot\cd\ot V_{q^{k(s-1)}},
\end{align*}
where $w_0$ is the longest element of $\mathfrak{S}_s$. Then $R_s$ is a $U'_q(\geh)$-linear
homomorphism. Define 
\[
V_s=\mbox{Im}\;R_s.
\]
Let us denote by $W$ the image of 
\[
R(q^k,q^{-k}):V_{q^k}\ot V_{q^{-k}}\longrightarrow V_{q^{-k}}\ot V_{q^k}
\]
and by $N$ its kernel. Then we have
\begin{align}
&V_s\text{ considered as a submodule of }V^{\ot s}=V_{q^{-k(s-1)}}\ot\cd\ot V_{q^{k(s-1)}}\\
&\text{is contained in }\bigcap_{i=0}^{s-2}V^{\ot i}\ot W\ot V^{\ot(s-2-i)}. \nonumber\\
\intertext{Similarly, we have}
&V_s\text{ is a quotient of }V^{\ot s}/\sum_{i=0}^{s-2}V^{\ot i}\ot N\ot V^{\ot(s-2-i)}.
\label{Vlquotient}
\end{align}

In the sequel, following \cite{KMN2} we define a prepolarization on $V_s$ and study
necessary properties. First we recall the following lemma.

\begin{lemma} \cite[Lemma 3.4.1]{KMN2}
Let $M_j$ and $N_j$ be $\Uqp$-modules and let $(\;,\;)_j$ be an admissible pairing between 
$M_j$ and $N_j$ $(j=1,2)$. Then the pairing $(\;,\;)$between $M_1\ot M_2$ and $N_1\ot N_2$
defined by $(u_1\ot u_2,v_1\ot v_2)=(u_1,v_1)_1(u_2,v_2)_2$ for all $u_j\in M_j$ and 
$v_j\in N_j$ is admissible.
\end{lemma}

Let $V$ be a finite-dimensional $\Uqp$-module satisfying \eqref{Virr} and \eqref{la0}. Suppose
$V$ has a polarization.
The polarization on $V$ gives an admissible pairing between $V_x$ and $V_{x^{-1}}$. Hence
it induces an admissible pairing between $V_{x_1}\ot\cd\ot V_{x_s}$ and 
$V_{x_1^{-1}}\ot\cd\ot V_{x_s^{-1}}$.

\begin{lemma} \cite[Lemma 3.4.2]{KMN2}
If $x_j=x_{s+1-j}^{-1}$ for $j=1,\ldots,s$, then for any $u,u'\in V_{x_1}\ot\cd\ot V_{x_s}$,
we have
\[
(u,R_{w_0}(x_1,\ldots,x_s)u')=(u',R_{w_0}(x_1,\ldots,x_s)u).
\]
\end{lemma}
By taking $x_i = q^{k(s-2i+1)}$, we obtain the admissible pairing $(\;,\;)$
between $W=V_{q^{k(s-1)}}\ot V_{q^{k(s-3)}}\ot\cd\ot V_{q^{-k(s-1)}}$ and 
$W'=V_{q^{-k(s-1)}}\ot V_{q^{-k(s-3)}}\ot\cd\ot V_{q^{k(s-1)}}$ that satisfies
\begin{equation} \label{wRw'}
(w,R_sw')=(w',R_sw)\quad\text{for any }w,w'\in W.
\end{equation}
This allows us to define a prepolarization $(\;,\;)_s$ on $V_s$ by
\[
(R_su,R_su')_s=(u,R_su')
\]
for $u,u'\in V_{q^{k(s-1)}}\ot V_{q^{k(s-3)}}\ot\cd\ot V_{q^{-k(s-1)}}$. 

Assume
\begin{equation} \label{Kzu0}
\text{$V$ admits a $\Uqp_{\Kz}$-submodule $V_{\Kz}$ such that $(V_{K_{\Z}})_{\la_0}=K_{\Z}u_0$.}
\end{equation}
Let us further set
\[
(V_s)_{K_{\Z}}=R_s((V_{K_{\Z}})^{\ot s})\cap(V_{K_{\Z}})^{\ot s}.
\]
Then \cite[Proposition3.4.3]{KMN2} follows:
\begin{proposition} \label{prop:3.3} \mbox{}
\begin{itemize}
\item[(i)] $(\;,\;)_s$ is a nondegenerate prepolarization on $V_s$.
\item[(ii)] $(R_s(u_0^{\ot s}),R_s(u_0^{\ot s}))_s=1$.
\item[(iii)] $((V_s)_{\Kz},(V_s)_{\Kz})_s\subset\Kz$.
\end{itemize}
\end{proposition}

\subsection{KR modules}

We want to apply the fusion construction with $V$ being the fundamental representation $W(\varpi_r)$.
Let us take $k$ to be $(\alpha_r,\alpha_r)/2$ except in the case $k=1$ when $\geh=A_{2n}^{(2)}$ and
$r=n$.

\begin{proposition} \label{prop:assump ok}
Assumptions \eqref{Virr},\eqref{la0},\eqref{pole} and \eqref{Kzu0} hold for the fundamental 
representations.
\end{proposition}

\begin{proof}
\eqref{Virr} is a consequence of Proposition 2.4 (v) and the fact that $B(\varpi_r)$ is a ``simple"
crystal (see \cite{K3}). \eqref{la0} is valid by Proposition 2.4 (iv) with $\la_0=cl(\varpi_r)$.
Noting that $W(\varpi_r)$ is a ``good" $\Uqp$-module, \eqref{pole} is the consequence of Proposition
9.3 of \cite{K3}. \eqref{Kzu0} is valid, since $W(\varpi_r)$ admits a $\Uqp_{\Z}$-submodule 
$W(\varpi_r)_{\Z}$ induced from $V(\varpi_r)_{\Z}$ such that $(W(\varpi_r)_{\Z})_{cl(\varpi_r)}
=\Z[q_s,q_s^{-1}]u_{\varpi_r}$.
\end{proof}

For $r\in I_0$ and $s\in\Z_{>0}$ we define the $\Uqp$-module $W^{(r)}_s$ to be the module
constructed by the fusion construction in section~\ref{ss:fusion} with $V=W(\varpi_r)$ and 
$k=(\alpha_r,\alpha_r)/2$ except in the case $k=1$ when $\geh=A_{2n}^{(2)}$ and $r=n$.

\begin{proposition} \label{prop:KRassump} \mbox{}
\begin{itemize}
\item[(i)] There exists a prepolarization $(\;,\;)$ on $W^{(r)}_s$.
\item[(ii)] There exists a $\Uqp_{\Kz}$-submodule $(W^{(r)}_s)_{\Kz}$ of $W^{(r)}_s$ such 
	that $$((W^{(r)}_s)_{\Kz},(W^{(r)}_s)_{\Kz})\subset\Kz.$$
\item[(iii)] There exists a vector $u_0$ of weight $s\varpi_r$ in $(W^{(r)}_s)_{\Kz}$ such that
	$(u_0,u_0)=1$.
\end{itemize}
\end{proposition}

\begin{proof}
The results follow from Propositions \ref{prop:3.3} and \ref{prop:assump ok}.
\end{proof}

The following proposition is an easy consequence of the main result of Kashiwara \cite{K3}.
Note also that his result can be applied not only to KR modules but also to any irreducible 
modules.

\begin{proposition} \label{prop:KR irred}
$W^{(r)}_s$ is irreducible and its Drinfeld polynomials are given by
\[
P_j(u)=\left\{
\begin{array}{ll}
(1-a_r^\dagger q_r^{1-s}u)(1-a_r^\dagger q_r^{3-s}u)\cd(1-a_r^\dagger q_r^{s-1}u)\quad&(j=r)\\
1&(j\neq r)
\end{array}\right.
\]
except when $\geh=A_{2n}^{(2)}$ and $r=n$. If $\geh=A_{2n}^{(2)}$ and $r=n$, they are given by 
replacing $q_r$ with $q$ in the above formula.
\end{proposition}

\begin{proof}
Let $V$ be a nonzero submodule of $V_s=W^{(r)}_s$. To show the irreducibility, it suffices to show 
that any vector $v$ in $V_s$ is 
contained in $V$. By definition there exists a vector $u\in W(\varpi_r)^{\ot s}$ such that $v=R_su$.
{}From Theorem 9.2 (ii) of \cite{K3} we have $u_0^{\ot s}\in V$. From Theorem 9.2 (i) of loc. cit. there
exists $x\in\Uqp$ such that $u=\Delta^{(s)}(x)u_0^{\ot s}$, where $\Delta^{(s)}$ is the coproduct
$\Uqp\longrightarrow\Uqp^{\ot s}$. Hence we have $v=R_s\Delta^{(s)}(x)u_0^{\ot s}=
\Delta^{(s)}(x)R_su_0^{\ot s}=\Delta^{(s)}(x)u_0^{\ot s}\in V$.

Since $W^{(r)}_s$ is the irreducible module in $(W_1^{(r)})_{q_r^{1-s}}\ot(W_1^{(s)})_{q_r^{3-s}}
\ot\cd\ot(W_1^{(r)})_{q_r^{s-1}}$ generated by $u_0^{\ot s}$, the latter statement is clear from 
\cite[Corollary 3.5]{CP}, Lemma \ref{lem:Dpoly} and the fact that if $V$ corresponds to $\{P_j(u)\}$,
then $V_a$ does to $\{P_j(au)\}$.
\end{proof}

This irreducible $\Uqp$-module $W^{(r)}_s$ is called Kirillov-Reshetikhin (KR) module.

Since the KR module $W^{(r)}_s$ is also a $U_q(\gc)$-module by restriction, we have the following
direct sum decomposition as a $U_q(\gc)$-module.
\begin{equation} \label{KRchar}
W^{(r)}_s\simeq\bigoplus_{\la\in\ol{P}_+}N_s^{(r)}(\la)\cdot\ol{V}(\la)
\end{equation}
Namely, $N^{(r)}_s(\la)$ is the multiplicity of the irreducible $U_q(\gc)$-module $\ol{V}(\la)$
in $W^{(r)}_s$. Then we have a criterion that the KR module has a crystal pseudobase.

\begin{proposition} \label{prop:main}
Suppose for any $\la\in\ol{P}_+$ such that $N_s^{(r)}(\la)>0$ there exist $u(\la)_j\in(W^{(r)}_s)_{\Kz}$
of weight $\la$ for $j=1,\ldots,N_s^{(r)}(\la)$. If we have $(u(\la)_j,u(\la)_k)\in\delta_{jk}+q_sA$ and
$(e_ju(\la)_k,e_ju(\la)_k)\in q_sq_j^{-2(1+\inner{h_j}{\la})}A$ for any $j\in I_0$, then $(\;,\;)$ on 
$W^{(r)}_s$ is a polarization, and $W^{(r)}_s$ has a crystal pseudobase.
\end{proposition}

\begin{proof}
We use Proposition \ref{prop:key}. All the assumptions except \eqref{norm} are satisfied by Propositions
\ref{prop:KRassump}. Note that $(u(\la)_j,u(\mu)_k)=0$ if $\la\neq\mu$.
\end{proof}

\begin{remark} \label{rem:existence}
{}From the previous proposition it immediately follows that if $W^{(r)}_s$ is irreducible as a 
$U_q(\geh_0)$-module, then it has a crystal pseudobase (see also \cite[Proposition 3.4.4]{KMN2}). 
There is another case in which the existence of crystal pseudobase is proven for any $l$ and any 
$\geh$ except $A_n^{(1)}$ as in~\cite[Proposition 3.4.5]{KMN2}. It corresponds to
$r=2$ when $\geh=B_n^{(1)}, D_n^{(1)}, A_{2n-1}^{(1)}$, $r=6$ when $\geh=E_6^{(1)}$,
and $r=1$ in all other cases. Here we follow the labeling of vertices of the Dynkin diagram by~\cite{Kac}. 
We remark that the crystal base of $W_1^{(r)}$ for such $r$ is treated in~\cite{BFKL}.
\end{remark}

There is an explicit formula of $N^{(r)}_s(\la)$ called the ($q=1$) fermionic formula. We have
\cite{Chari,HKOTT,HKOTY,H1,H2,KuNa,N2,N3} for references. 
To explain it, we introduce $t_i$ and $t_i^\vee$ for $i\in I_0$ by 
\[
t_i=\left\{
\begin{array}{ll}
\frac{2}{(\alpha_i,\alpha_i)}\quad&\text{if $\geh$ is untwisted}\\
1&\text{if $\geh$ is twisted}
\end{array}\right.
\]
and $t_i^\vee=(t_i\text{ for }\geh^\vee)$, where $\geh^\vee$ is the dual Kac-Moody algebra to $\geh$. For 
$p\in\Z$ and $m\in\Z_{\ge0}$ let ${p+m\choose m}$ stand for the binomial coefficient, i.e., 
${p+m\choose m}=\prod_{k=1}^m\frac{p+k}{k}$. Then, for $r\in I_0,s\in\Z_{>0}$ and $\la\in\ol{P}_+$ 
we have
\begin{align*}
N_s^{(r)}(\la)&=\sum_{\bf m}\prod_{a\in I_0,j\ge1}{p_j^{(a)}+m_j^{(a)}\choose m_j^{(a)}},\\
\intertext{where}
p_j^{(a)}&=\delta_{ai}\min(j,s)-\frac1{t_a^\vee}\sum_{b\in I_0,k\ge1}
(\alpha_a,\alpha_b)\min(t_bj,t_ak)m_k^{(b)}
\end{align*}
and the sum $\sum_{\bf m}$ is taken over all $(m_j^{(a)}\in\Z_{\ge0}\mid a\in I_0,j\ge1)$ satisfying
\[
\sum_{a\in I_0,j\ge1}jm_j^{(a)}\alpha_a=s\varpi_r-\la.
\]

The proof of this formula goes as follows.
Set $Q^{(r)}_s=\mathrm{ch}\,W^{(r)}_s$. It suffices to show that 
$Q^{(r)}_s=\sum_{\la\in\ol{P}_+}N_s^{(r)}(\la)\mathrm{ch}\,\ol{V}(\la)$.
By Theorem 8.1 of \cite{HKOTY} (see also Theorem 6.3 of \cite{HKOTT} including the twisted
cases), it suffices to show that $\{Q^{(r)}_s\}$ satisfies the conditions (A),(B),(C) in 
the theorem. (A) is evident by the construction of $W^{(r)}_s$, and (B),(C) were verified 
in \cite{N2,H1,H2} for the simply-laced, untwisted and twisted cases, respectively. Note that
condition (C) is replaced with another convergence property (4.15) of \cite{KNT}.
Note also that there is an earlier result by Chari \cite{Chari} for untwisted cases.
It should also be noted that there is another explicit formula
$M_s^{(r)}(\lambda)$ for the multiplicities $N_s^{(r)}(\lambda)$ which involves unsigned
binomial coefficients, that is ${p+m\choose m}=0$ if $p<0$~\cite{HKOTY,HKOTT}.
It was recently shown by Di Francesco and Kedem~\cite{DK:2007} that $M_s^{(r)}(\lambda)
=N_s^{(r)}(\lambda)$ in the untwisted cases.

For nonexceptional types, the explicit value of $N^{(r)}_s(\lambda)$ can be found in section
7 of \cite{HKOTY} for untwisted cases, and in section 6.2 of \cite{HKOTT} for twisted cases.
See \eqref{branching}.

\section{Existence of crystal pseudobases for nonexceptional types}
\label{sec:existence}

In this section we show that any KR module for nonexceptional type has a crystal pseudobase.
For type $A_n^{(1)}$ this fact is established in \cite{KMN2}. So we do not deal with the $A_n^{(1)}$
case.

\subsection{Dynkin data} \label{subsec:Dynkin}

\begin{table}
\begin{align*}
&\xymatrix@R=1ex{
& *{\circ}<3pt> \ar@{-}[dr]^<{0} \\
B_n^{(1)} & & *{\circ}<3pt> \ar@{-}[r]_<{2} 
& {} \ar@{.}[r]&{}  \ar@{-}[r]_>{\,\,\,\,n-1} &
*{\circ}<3pt> \ar@{=}[r] |-{\scalebox{2}{\object@{>}}}& *{\bullet}<3pt>\ar@{}_<{n} & (\vdom,B_n)\\
& *{\circ}<3pt> \ar@{-}[ur]_<{1}
}\\
&\xymatrix{
C_n^{(1)}& *{\circ}<3pt> \ar@{=}[r] |-{\scalebox{2}{\object@{>}}}_<{0} 
&*{\circ}<3pt> \ar@{-}[r]_<{1} 
& {} \ar@{.}[r]&{}  \ar@{-}[r]_>{\,\,\,\,n-1} &
*{\circ}<3pt> \ar@{=}[r] |-{\scalebox{2}{\object@{<}}}
& *{\bullet}<3pt>\ar@{}_<{n} & (\hdom,C_n)
}\\
&\xymatrix@R=1ex{
& *{\circ}<3pt> \ar@{-}[dr]^<{0}&&&&&*{\bullet}<3pt> \ar@{-}[dl]^<{n-1}\\
D_n^{(1)} & & *{\circ}<3pt> \ar@{-}[r]_<{2} 
& {} \ar@{.}[r]&{} \ar@{-}[r]_>{\,\,\,n-2} &
*{\circ}<3pt> && (\vdom,D_n)\\
& *{\circ}<3pt> \ar@{-}[ur]_<{1}&&&&&
*{\bullet}<3pt> \ar@{-}[ul]_<{n}
}\\
&\xymatrix{
A_{2n}^{(2)}& *{\circ}<3pt> \ar@{=}[r] |-{\scalebox{2}{\object@{<}}}_<{0} 
&*{\circ}<3pt> \ar@{-}[r]_<{1} 
& {} \ar@{.}[r]&{}  \ar@{-}[r]_>{\,\,\,\,n-1} &
*{\circ}<3pt> \ar@{=}[r] |-{\scalebox{2}{\object@{<}}}
& *{\circ}<3pt>\ar@{}_<{n} & (\cell,C_n)
}\\
&\xymatrix@R=1ex{
& *{\circ}<3pt> \ar@{-}[dr]^<{0} \\
A_{2n-1}^{(2)}& & *{\circ}<3pt> \ar@{-}[r]_<{2} 
& {} \ar@{.}[r]&{}  \ar@{-}[r]_>{\,\,\,\,n-1} &
*{\circ}<3pt> \ar@{=}[r] |-{\scalebox{2}{\object@{<}}}& *{\circ}<3pt>\ar@{}_<{n} & (\vdom,C_n)\\
& *{\circ}<3pt> \ar@{-}[ur]_<{1}
}\\
&\xymatrix{
D_{n+1}^{(2)}& *{\circ}<3pt> \ar@{=}[r] |-{\scalebox{2}{\object@{<}}}_<{0} 
&*{\circ}<3pt> \ar@{-}[r]_<{1} 
& {} \ar@{.}[r]&{}  \ar@{-}[r]_>{\,\,\,\,n-1} &
*{\circ}<3pt> \ar@{=}[r] |-{\scalebox{2}{\object@{>}}}
& *{\bullet}<3pt>\ar@{}_<{n} & (\cell,B_n)
}
\end{align*}
\caption{\label{tab:Dynkin}Dynkin diagrams}
\end{table}
First we list the Dynkin diagrams of all nonexceptional affine algebras except $A_n^{(1)}$ in Table
\ref{tab:Dynkin}. We also list the pair $(\nu,\geh_0)$ in the table with a partition $\nu=\vdom,\hdom,\cell$
and a simple Lie algebra $\geh_0$ whose Dynkin diagram is the one obtained by removing the $0$-vertex.
Note that the difference of $\nu$ comes from the diagram near the $0$-vertex.

The simple roots for type $B_n,C_n,D_n$ are
\begin{equation*}
\begin{split}
\alpha_i&=\epsilon_i-\epsilon_{i+1}  \qquad \text{for $1\le i<n$}\\
\alpha_n&= \begin{cases} 
		\epsilon_{n-1}+\epsilon_n & \text{for type $D_n$}\\
		\epsilon_n                             &\text{for type $B_n$}\\
		2 \epsilon_n                          & \text{for type $C_n$}
		\end{cases}
\end{split}
\end{equation*}
and the fundamental weights are
\begin{equation*}
\begin{aligned}[3]
\text{Type $D_n$:} \qquad
	&\varpi_i = \epsilon_1+\cdots+\epsilon_i &&\text{for $1\le i\le n-2$}\\
	&\varpi_{n-1} = (\epsilon_1+\cdots+\epsilon_{n-1}-\epsilon_n)/2 &&\\
	&\varpi_n = (\epsilon_1+\cdots+\epsilon_{n-1}+\epsilon_n)/2&&\\[2mm]
\text{Type $B_n$:} \qquad
	&\varpi_i = \epsilon_1+\cdots+\epsilon_i &&\text{for $1\le i\le n-1$}\\
	&\varpi_n = (\epsilon_1+\cdots+\epsilon_{n-1}+\epsilon_n)/2&&\\[2mm]
\text{Type $C_n$:} \qquad
	&\varpi_i = \epsilon_1+\cdots+\epsilon_i &&\text{for $1\le i\le n$}
\end{aligned}
\end{equation*}
where $\epsilon_i$ ($i=1,\ldots,n$) are vectors in the weight space of each simple Lie algebra.
(By convention we set $\varpi_0=0$.)
These elements can be viewed as those of the weight lattice $P$ of the affine algebra in Table 
\ref{tab:Dynkin}. On $P$ we defined the inner product $(\;,\;)$ normalized as $(\delta,\la)=
\inner{c}{\la}$ for $\la\in P$. This normalization is equivalent to setting $(\epsilon_i,\epsilon_j)
=\kappa\delta_{ij}$ with $\kappa=\frac12$ for $C_n^{(1)}$, $=2$ for $D_{n+1}^{(2)}$, and $=1$ for the
other types. However, in this section we renormalize it by $(\epsilon_i,\epsilon_j)=\delta_{ij}$.
This is equivalent to setting $(\alpha_i,\alpha_i)/2=1$ for $i$ not an end node of the Dynkin diagram.
We also note that 
\[
\alpha_0=\left\{
\begin{array}{ll}
\delta-\epsilon_1-\epsilon_2\quad&\mbox{if }\nu=\vdom\\
\delta-2\epsilon_1&\mbox{if }\nu=\hdom\\
\delta-\epsilon_1&\mbox{if }\nu=\cell\;.
\end{array}
\right.
\]

\subsection{Existence of crystal pseudobases for KR modules}

We first present the branching rule of KR modules of affine type listed in Table \ref{tab:Dynkin}
with respect to the subalgebra $U_q(\geh_0)$. They can be found 
in~\cite[Theorems 7.1 and 8.1]{HKOTY} and~\cite[Theorems 6.2 and 6.3]{HKOTT}.
For $i\in I_0$ for $\geh$ we say $i$ is a spin node
if the vertex $i$ is filled in Table \ref{tab:Dynkin}. If $r\in I_0$ is a spin node, then the KR module
$W_s^{(r)}$ is irreducible as a $U_q(\geh_0)$-module:
\[
W_s^{(r)}\simeq\ol{V}(s\varpi_r).
\]
Suppose now that $r\in I_0$ is not a spin node. Let $\om$ be a dominant integral weight of the form
of $\om=\sum_ic_i\varpi_i$. Assume $c_i=0$ for $i$ a spin node. In the standard way we represent $\om$
by the partition that has exactly $c_i$ columns of height $i$. Then the KR module $W_s^{(r)}$ 
decomposes into
\begin{equation} \label{branching}
W_s^{(r)}\simeq\bigoplus_\om\ol{V}(\om)
\end{equation}
as a $U_q(\geh_0)$-module, where $\om$ runs over all partitions that can be obtained from the $r\times s$
rectangle by removing pieces of shape $\nu$ (with $\nu$ as in Table~\ref{tab:Dynkin}).

If $r\in I_0$ is a spin node, the KR module $W_s^{(r)}$ has a crystal pseudobase by Remark 
\ref{rem:existence}. Suppose $r$ is not a spin node. As we have seen, we have $N_s^{(r)}(\la)\le1$.
Hence, by Proposition~\ref{prop:main}, in order to show the existence of crystal pseudobase, it 
suffices to define a vector 
$u(\la)\in(W_s^{(r)})_{\Kz}$ of weight $\la$ for any $\la$ such that $N_s^{(r)}=1$, and show 
$(u(\la),u(\la))\in 1+q_sA$ and $(e_ju(\la),e_ju(\la))\in q_sq_j^{-2(1+\inner{h_j}{\la})}A$ for $j\in I_0$.
In the subsequent subsections, we do this task by dividing into 3 cases according to the shape of $\nu$.

\subsection{Calculation of prepolarization: $D_n^{(1)},B_n^{(1)},A_{2n-1}^{(2)}$ cases}

We assume $1\le r\le n-2$ for $D_n^{(1)}$, $1\le r\le n-1$ for $B_n^{(1)}$
and $1\le r\le n$ for $A_{2n-1}^{(2)}$. Let $r'=[r/2]$. 
Let $\cb=(c_1,c_2,\ldots,c_{r'})$ be a sequence of integers such that 
$s\ge c_1\ge c_2\ge\cd\ge c_{r'}\ge0$. For such $\cb$ we define a vector $u_m$ ($0\le m\le r'$) 
in $W_s^{(r)}$ inductively by
\[
u_m=(e_{r-2m}^{(c_m)}\cd e_2^{(c_m)}e_1^{(c_m)})(e_{r-2m+1}^{(c_m)}\cd e_3^{(c_m)}e_2^{(c_m)})
e_0^{(c_m)}u_{m-1},
\]
where $u_0$ is the vector in (iii) of Proposition \ref{prop:KRassump}. Set $u(\cb)=u_{r'}$. The weight of 
$u(\cb)$ is given by
\[
\la(\cb)=\sum_{j=0}^{r'}(c_j-c_{j+1})\varpi_{r-2j},
\]
where we have set $c_0=s,c_{r'+1}=0$, and $\varpi_0$ should be understood as $0$. $\la(\cb)$ represents
all $\om$ in \eqref{branching} when $\cb$ runs over all possible sequences.
For $l,m\in\Z_{\ge0}$ such that $m\le l$ we define the $q$-binomial coefficient by
\begin{equation} \label{q-binom}
{l\brack m}=\frac{[l]!}{[m]![l-m]!}.
\end{equation}
The following proposition calculates values of the prepolarization $(\;,\;)$ on $W_s^{(r)}$.
\begin{proposition} \label{prop:norm DBA} \mbox{}
\begin{itemize}
\item[(1)] ${\displaystyle (u(\cb),u(\cb))=\prod_{j=1}^{r'}q^{c_j(2s-c_j)}{2s\brack c_j}}$,
\item[(2)] $(e_ju(\cb),e_ju(\cb))=0$ unless $r-j\in 2\Z_{\ge0}$. If $r-j\in 2\Z_{\ge0}$, then 
	setting $p=(r-j)/2+1$, $(e_ju(\cb),e_ju(\cb))$ is given by
\[
q^{2s-c_{p-1}-1}[2s-c_{p-1}]\prod_{j=1}^{r'}q^{(c_j-\delta_{j,p})(2s-c_j)}
{2s-\delta_{j,p}\brack c_j-\delta_{j,p}}.
\]
\end{itemize}
\end{proposition}

For type $D_n^{(1)}$ this proposition is proven in~\cite{O:2006}. The proof goes completely 
parallel also for type $B_n^{(1)}$ and $A_{2n-1}^{(2)}$. Note that $q_i=q$ for $i\neq n$, 
$q_n=q,q^{1/2},q^2$ for $D_n^{(1)},B_n^{(1)},A_{2n-1}^{(2)}$, respectively, and $q_s=q^{1/2}$ for 
$B_n^{(1)}$, $=q$ for $D_n^{(1)},A_{2n-1}^{(2)}$. Since $q^{m-1}[m],q^{n(m-n)}\in1+qA$ and 
$\inner{h_j}{\la(\cb)}=c_{p-1}-c_p\ge0$, we have $(u(\cb),u(\cb))\in1+q_sA$ and $(e_ju(\cb),e_ju(\cb))
\in q_sq_j^{-2(1+\inner{h_j}{\la(\cb)})}A$. for $j\in I_0$. This establishes the conditions of
Proposition~\ref{prop:main} and hence proves Theorem~\ref{thm:main} that $W^{(r)}_s$ has a 
crystal pseudobase.

We denote the crystal of $W_s^{(r)}$ by $B^{r,s}$.
Similar to $\geh_0$ one can consider $\geh_1$, which is another (mutually isomorphic) simple Lie
algebra obtained by removing the vertex $1$ from the Dynkin diagram of $\geh$.
The following proposition will be used to show that $B^{r,s}$ is isomorphic to $\Bt^{r,s}$, which 
is given combinatorially in the next section.

\begin{proposition} \label{prop:classical decomp}
Let $1\le r\le n-2$ for $\geh=D_n^{(1)}$, $1\le r\le n-1$ for $\geh=B_n^{(1)}$, 
$1\le r\le n$ for $\geh=A_{2n-1}^{(2)}$, and $s\in\Z_{>0}$.
Then for $i=0,1$, $B^{r,s}$ decomposes as $U_q(\geh_i)$-crystals into 
\[
B^{r,s}\simeq\bigoplus_{0\le m_1\le\cd\le m_s\le [r/2]}
B^{\geh_i}(\sigma^i(\varpi_{r-2m_1}+\cd+\varpi_{r-2m_s})).
\]
Here $B^{\geh_i}(\la)$ is the crystal base of the highest weight $U_q(\geh_i)$-module of highest weight 
$\la$, and $\sigma$ is the automorphism on $P$ such that $\sigma(\La_0)=\La_1,\sigma(\La_1)=\La_0,
\sigma(\La_j)=\La_j$ ($j>1$) and extended linearly.
\end{proposition}

\begin{proof}
If $i=0$, the claim is a direct consequence of \eqref{branching}. For $i=1$ note that the Weyl group
of $\gc$ contains an element $w$ which sends $\varpi_j$ to $\sigma(\varpi_j)$ for any $j$ such that
$0\le j\le r$, where by convention $\varpi_0=0$. (Using the orthogonal basis $\{\epsilon_i\}$ of 
section \ref{subsec:Dynkin} of the weight space
of $\gc$, we can take an element $w$ such that $w(\epsilon_i)=(-1)^{\delta(i)}\epsilon_i$, where
$\delta(i)=1$ if $i=1,n$ for $\geh=D_n^{(1)}$, $i=1$ for $\geh=B_n^{(1)}$ and $A_{2n-1}^{(2)}$, and
$\delta(i)=0$ otherwise.) Since $W^{(r)}_s$ is a direct sum also as a $U_q(\geh_1)$-module, it is 
enough to show the following equality of characters.
\begin{equation} \label{proof}
\mathrm{ch}\,W^{(r)}_s=\sum_{0\le m_1\le\cd\le m_s\le [r/2]}
\mathrm{ch}\,V^{\geh_1}(\sigma(\varpi_{r-2m_1}+\cd+\varpi_{r-2m_s}))
\end{equation}
Here $V^{\geh_1}(\la)$ denotes the highest weight $U_q(\geh_1)$-module of highest weight $\la$. 
But noting $w(\alpha_0)=\alpha_1,w(\alpha_1)=\alpha_0,w(\alpha_j)=\alpha_j$ ($j>1$) on $P_{cl}$,
\eqref{proof} is shown from 
\[
\mathrm{ch}\,W^{(r)}_s=\sum_{0\le m_1\le\cd\le m_s\le [r/2]}
\mathrm{ch}\,V^{\geh_0}(\varpi_{r-2m_1}+\cd+\varpi_{r-2m_s})
\]
since $w$ preserves the weight multiplicity.
\end{proof}

\subsection{Calculation of prepolarization: $C_n^{(1)}$ case}

We assume $1\le r\le n-1$.
Let $\cb=(c_1,c_2,\ldots,c_r)$ be a sequence of integers such that 
$[s/2]\ge c_1\ge c_2\ge\cd\ge c_r\ge0$. For such $\cb$ we define a vector $u_m$ ($0\le m\le r$) 
in $W_s^{(r)}$ inductively by
\[
u_m=e_{r-m}^{(2c_m)}\cd e_2^{(2c_m)}e_1^{(2c_m)}e_0^{(c_m)}u_{m-1},
\]
where $u_0$ is the vector in (iii) of Proposition \ref{prop:KRassump}. Set $u(\cb)=u_r$. The weight of 
$u(\cb)$ is given by
\[
\la(\cb)=\sum_{j=0}^r2(c_j-c_{j+1})\varpi_{r-j},
\]
where we have set $c_0=s/2,c_{r+1}=0$, and $\varpi_0$ should be understood as $0$. $\la(\cb)$ represents
all $\om$ in \eqref{branching} when $\cb$ runs over all possible sequences. In this subsection, besides
\eqref{q-binom} we also use ${l\brack m}_0$ defined by \eqref{q-binom} with $q$ replaced by $q_0=q^2$.
(Recall that we have renormalized the inner product $(\;,\;)$ on $P$ in such a way that 
$(\epsilon_i,\epsilon_j)=\delta_{ij}$.)

We are to calculate the values of $(u(\cb),u(\cb))$ and $(e_ju(\cb),e_ju(\cb))$. Since the calculation
goes parallel to the case of $D_n^{(1)}$ treated in \cite{O:2006}, we only give here intermediate results 
as a lemma. We write $\Vert u\Vert^2$ for $(u,u)$.

\begin{lemma} \mbox{}
\begin{itemize}
\item[(1)] $\Vert u_m\Vert^2=q_0^{c_m(s-c_m)}{s\brack c_m}_0\Vert u_{m-1}\Vert^2$,
\item[(2)] $e_ju(\cb)=0$ if $j>r$,
\item[(3)] $\Vert e_ju(\cb)\Vert^2=q^{2\beta_j}\Vert f_ju(\cb)\Vert^2+q^{\beta_j-1}[\beta_j]
	\Vert u(\cb)\Vert^2$ if $1\le j\le r$, where $\beta_j=-\inner{h_j}{\la(\cb)}=2(c_{r+1-j}-c_{r-j})$,
\item[(4)] 
	\begin{equation*}
	\begin{split}
	\Vert f_ju(\cb)\Vert^2=& \prod_{\substack{1\le m\le r\\ m\neq r-j+1}}
	q_0^{c_m(s-c_m)}{s\brack c_m}_0\\
	&\times q_0^{c_{r-j+1}(s-1-c_{r-j+1})}{s-1\brack c_{r-j+1}}_0\times q^{2c_{r-j}-1}[2c_{r-j}].
	\end{split}
	\end{equation*}
\end{itemize}
\end{lemma}
{}From this lemma we have

\begin{proposition} \label{prop:norm C}\mbox{}
\begin{itemize}
\item[(1)] $(u(\cb),u(\cb))=\prod_{m=1}^r q^{c_m(s-c_m)}{s\brack c_m}_0$,
\item[(2)] 
\[(e_ju(\cb),e_ju(\cb))= \begin{cases}
	q^{2s-2c_{r-j}-1}[2s-2c_{r-j}] &\\
	\quad \times \prod_{m=1}^rq_0^{(c_m-\delta_{m,r-j+1})(s-c_m)}
	{s-\delta_{m,r-j+1}\brack c_m-\delta_{m,r-j+1}} & \text{if $1\le j\le r$}\\
	0& \text{if $r<j\le n$}.
	\end{cases}
\]
\end{itemize}
\end{proposition}
Note that $q_i=q$ for $i\neq0,n$, $q_n=q^2$, and $q_s=q$ under the renormalization. Since 
$\inner{h_j}{\la(\cb)}=-\beta_j=2(c_{r-j}-c_{r+1-j})\ge0$, we have $(u(\cb),u(\cb))\in1+q_sA$ and
$(e_ju(\cb),e_ju(\cb))\in q_sq_j^{-2(1+\inner{h_j}{\la(\cb)})}A$ for $j\in I_0$. 
By Proposition~\ref{prop:main} this proves Theorem~\ref{thm:main}.

\subsection{Calculation of prepolarization: $A_{2n}^{(2)},D_{n+1}^{(2)}$ cases}

We assume $1\le r\le n$ for $A_{2n}^{(2)}$ and $1\le r\le n-1$ for $D_{n+1}^{(2)}$.
Let $\cb=(c_1,c_2,\ldots,c_r)$ be a sequence of integers such that 
$s\ge c_1\ge c_2\ge\cd\ge c_r\ge0$. For such $\cb$ we define a vector $u_m$ ($0\le m\le r$) 
in $W_s^{(r)}$ inductively by
\[
u_m=e_{r-m}^{(c_m)}\cd e_1^{(c_m)}e_0^{(c_m)}u_{m-1},
\]
where $u_0$ is the vector in (iii) of Proposition \ref{prop:KRassump}. Set $u(\cb)=u_r$. The weight of 
$u(\cb)$ is given by
\[
\la(\cb)=\sum_{j=0}^r(c_j-c_{j+1})\varpi_{r-j},
\]
where we have set $c_0=s,c_{r+1}=0$, and $\varpi_0$ should be understood as $0$. $\la(\cb)$ represents
all $\om$ in \eqref{branching} when $\cb$ runs over all possible sequences. In this subsection, besides
\eqref{q-binom} we also use ${l\brack m}_0$ defined by \eqref{q-binom} with $q$ replaced by $q_0=q^{1/2}$.

As in the previous subsection, we only give here intermediate results 
as a lemma. As before we write $\Vert u\Vert^2$ for $(u,u)$.

\begin{lemma} \mbox{}
\begin{itemize}
\item[(1)] $\Vert u_m\Vert^2=q_0^{c_m(2s-c_m)}{2s\brack c_m}_0\Vert u_{m-1}\Vert^2$,
\item[(2)] $e_ju(\cb)=0$ if $j>r$,
\item[(3)] $\Vert e_ju(\cb)\Vert^2=q^{2\beta_j}\Vert f_ju(\cb)\Vert^2+q^{\beta_j-1}[\beta_j]
	\Vert u(\cb)\Vert^2$ if $1\le j\le r$, where $\beta_j=-\inner{h_j}{\la(\cb)}=c_{r+1-j}-c_{r-j}$,
\item[(4)]
\begin{align*}
\Vert f_ju(\cb)\Vert^2&=\prod_{m=1}^rq_0^{c_m(2s-2\delta^{(1)}-c_m)}{2s-2\delta^{(1)}\brack c_m}_0
\times q^{c_{r-j}-1}[c_{r-j}]\\
&+\prod_{m=1}^rq_0^{(c_m+\delta^{(1)}-\delta^{(2)})(2s-\delta^{(1)}+\delta^{(2)}-c_m)}
{2s-2\delta^{(1)}\brack c_m-\delta^{(1)}-\delta^{(2)}}_0\times [2s-c_{r-j}+1]_0^2,
\end{align*}
where $\delta^{(1)}=\delta_{m,r-j+1},\delta^{(2)}=\delta_{m,r-j}$.
\end{itemize}
\end{lemma}
{}From this lemma we have

\begin{proposition} \label{prop:norm AD}\mbox{}
\begin{itemize}
\item[(1)] $(u(\cb),u(\cb))=\prod_{m=1}^r q^{c_m(2s-c_m)}{2s\brack c_m}_0$,
\item[(2)] 
\[(e_ju(\cb),e_ju(\cb))=\begin{cases}
	q^{2\beta_j}\Vert f_ju(\cb)\Vert^2+q^{\beta_j-1}[\beta_j]\Vert u(\cb)\Vert^2
	& \text{if $1\le j\le r$}\\
	0 & \text{if $r<j\le n$,}
	\end{cases}
\]
where $\beta_j$ and $\Vert f_ju(\cb)\Vert^2$ are given in the previous lemma.
\end{itemize}
\end{proposition}
Note that $q_i=q$ for $i\neq0,n$, $q_n=q^2$ for $A_{2n}^{(2)}$, $=q^{1/2}$ for $D_{n+1}^{(2)}$, 
and $q_s=q^{1/2}$ under the renormalization. Since 
$\inner{h_j}{\la(\cb)}=-\beta_j=c_{r-j}-c_{r+1-j}\ge0$, we have $(u(\cb),u(\cb))\in1+q_sA$ and
$(e_ju(\cb),e_ju(\cb))\in q_sq_j^{-2(1+\inner{h_j}{\la(\cb)})}A$ for $j\in I_0$.
By Proposition~\ref{prop:main} this proves Theorem~\ref{thm:main}.

\section{Combinatorial crystal $\Bt^{r,s}$ of type $D_n^{(1)},B_n^{(1)},A_{2n-1}^{(2)}$}
\label{sec:combinatorial KR}

In this section we review the combinatorial crystal $\Bt^{r,s}$ of~\cite{S:2007,St:2006} of type
$D_n^{(1)}$, $B_n^{(1)}$, and $A_{2n-1}^{(2)}$ and prove
some preliminary results that will be needed in section~\ref{sec:equivalence} to establish the
equivalence of $\Bt^{r,s}$ and $B^{r,s}$.

\subsection{Type $D_n$, $B_n$, and $C_n$ crystals} \label{ss:classical}

Crystals associated with a $U_q(\geh)$-module when $\geh$ is a simple Lie algebra of 
nonexceptional type, were studied by Kashiwara and Nakashima~\cite{KN:1994}.
Here we review the combinatorial structure in terms of tableaux of the crystals of type $X_n=
D_n$, $B_n$, and $C_n$ since these are the finite subalgebras relevant to the KR crystals of 
type $D_n^{(1)}$,  $B_n^{(1)}$, and $A_{2n-1}^{(2)}$.

For $\geh=D_n^{(1)},B_n^{(1)}$, or $A_{2n-1}^{(2)}$, any $\geh_0$ dominant weight $\om$ without 
a spin component can be expressed as $\om=\sum_i c_i \varpi_i$ for nonnegative integers $c_i$ 
and the sum runs over all $i=1,2,\ldots,n$ not a spin node. As explained earlier we represent $\om$ 
by the partition that has exactly $c_i$ columns of height $i$. 
For type $D_n$, this can be extended by associating a column of height $n-1$ with 
$\varpi_{n-1}+\varpi_n$ and a column of height $n$ with $2\varpi_n$. For type $B_n$
one may associate a column of height $n$ with $2\varpi_n$.
Conversely, if $\om$ is a partition, we write $c_i(\om)$ for the number of columns of $\om$ of 
height $i$. From now on we identify partitions and dominant weights in this way.

The crystal graph $B(\varpi_1)$ of the vector representation for type $D_n$, $B_n$, and
$C_n$ is given in Table \ref{tab:vr} by removing the 0 arrows in the crystal $B^{1,1}$
of type $D_n^{(1)}$, $B_n^{(1)}$, and $A_{2n-1}^{(2)}$, respectively.
\begin{table}
\begin{tabular}{|c|l|}
\hline
%
$D_n^{(1)}$ & \raisebox{-1.3cm}{\scalebox{0.7}{
\begin{picture}(365,100)(-10,-50)
\BText(0,0){1} \LongArrow(10,0)(30,0) \BText(40,0){2}
\LongArrow(50,0)(70,0) \Text(85,0)[]{$\cdots$}
\LongArrow(95,0)(115,0) \BText(130,0){n-1}
\LongArrow(143,2)(160,14) \LongArrow(143,-2)(160,-14)
\BText(170,15){n} \BBoxc(170,-15)(13,13)
\Text(170,-15)[]{\footnotesize$\overline{\mbox{n}}$}
\LongArrow(180,14)(197,2) \LongArrow(180,-14)(197,-2)
\BBoxc(215,0)(25,13)
\Text(215,0)[]{\footnotesize$\overline{\mbox{n-1}}$}
\LongArrow(230,0)(250,0) \Text(265,0)[]{$\cdots$}
\LongArrow(275,0)(295,0) \BBoxc(305,0)(13,13)
\Text(305,0)[]{\footnotesize$\overline{\mbox{2}}$}
\LongArrow(315,0)(335,0) \BBoxc(345,0)(13,13)
\Text(345,0)[]{\footnotesize$\overline{\mbox{1}}$}
\LongArrowArc(192.5,-367)(402,69,111)
\LongArrowArcn(152.5,367)(402,-69,-111) \PText(20,2)(0)[b]{1}
\PText(60,2)(0)[b]{2} \PText(105,2)(0)[b]{n-2}
\PText(152,13)(0)[br]{n-1} \PText(152,-9)(0)[tr]{n}
\PText(188,13)(0)[bl]{n} \PText(188,-9)(0)[tl]{n-1}
\PText(240,2)(0)[b]{n-2} \PText(285,2)(0)[b]{2}
\PText(325,2)(0)[b]{1} \PText(192.5,38)(0)[b]{0}
\PText(152.5,-35)(0)[t]{0}
\end{picture}
}}
\\ \hline
%
$B_n^{(1)}$ & \raisebox{-1.3cm}{\scalebox{0.7}{
\begin{picture}(350,100)(-10,-50)
\BText(0,0){1} \LongArrow(10,0)(30,0) \BText(40,0){2}
\LongArrow(50,0)(70,0) \Text(85,0)[]{$\cdots$}
\LongArrow(95,0)(115,0) \BText(125,0){n} \LongArrow(135,0)(155,0)
\BText(165,0){0} \LongArrow(175,0)(195,0) \BBoxc(205,0)(13,13)
\Text(205,0)[]{\footnotesize$\overline{\mbox{n}}$}
\LongArrow(215,0)(235,0) \Text(250,0)[]{$\cdots$}
\LongArrow(260,0)(280,0) \BBoxc(290,0)(13,13)
\Text(290,0)[]{\footnotesize$\overline{\mbox{2}}$}
\LongArrow(300,0)(320,0) \BBoxc(330,0)(13,13)
\Text(330,0)[]{\footnotesize$\overline{\mbox{1}}$}
\LongArrowArc(185,-330)(365,68,112)
\LongArrowArcn(145,330)(365,-68,-112) \PText(20,2)(0)[b]{1}
\PText(60,2)(0)[b]{2} \PText(105,2)(0)[b]{n-1}
\PText(145,2)(0)[b]{n} \PText(185,2)(0)[b]{n}
\PText(225,2)(0)[b]{n-1} \PText(270,2)(0)[b]{2}
\PText(310,2)(0)[b]{1} \PText(185,38)(0)[b]{0}
\PText(145,-35)(0)[t]{0}
\end{picture}
}}
\\ \hline
%
$A_{2n-1}^{(2)}$ & \raisebox{-1.3cm}{\scalebox{0.7}{
\begin{picture}(310,100)(-10,-50)
\BText(0,0){1} \LongArrow(10,0)(30,0) \BText(40,0){2}
\LongArrow(50,0)(70,0) \Text(85,0)[]{$\cdots$}
\LongArrow(95,0)(115,0) \BText(125,0){n} \LongArrow(135,0)(155,0)
\BBoxc(165,0)(13,13)
\Text(165,0)[]{\footnotesize$\overline{\mbox{n}}$}
\LongArrow(175,0)(195,0) \Text(210,0)[]{$\cdots$}
\LongArrow(220,0)(240,0) \BBoxc(250,0)(13,13)
\Text(250,0)[]{\footnotesize$\overline{\mbox{2}}$}
\LongArrow(260,0)(280,0) \BBoxc(290,0)(13,13)
\Text(290,0)[]{\footnotesize$\overline{\mbox{1}}$}
\LongArrowArc(165,-240)(275,65,115)
\LongArrowArcn(125,240)(275,-65,-115) \PText(20,2)(0)[b]{1}
\PText(60,2)(0)[b]{2} \PText(105,2)(0)[b]{n-1}
\PText(145,2)(0)[b]{n} \PText(185,2)(0)[b]{n-1}
\PText(230,2)(0)[b]{2} \PText(270,2)(0)[b]{1}
\PText(165,38)(0)[b]{0} \PText(125,-35)(0)[t]{0}
\end{picture}
}}
\\ \hline
\end{tabular}\vspace{4mm}
\caption{\label{tab:vr}KR crystal $B^{1,1}$}
\end{table}
The crystal $B(\varpi_\ell)$ for $\ell$ not a spin node can be realized as the connected component
of $B(\varpi_1)^{\otimes \ell}$ containing the element $\ell\otimes (\ell-1)\otimes \cdots \otimes 1$,
where we use the anti-Kashiwara convention for tensor products.
Similarly, the crystal $B(\om)$ labeled by a dominant weight $\om=\varpi_{\ell_1}+\cdots+\varpi_{\ell_k}$
with $\ell_1\ge \ell_2 \ge \cdots \ge \ell_k$ not containing spin nodes can be realized as the connected 
component in $B(\varpi_{\ell_1}) \otimes \cdots \otimes B(\varpi_{\ell_k})$ containing the element 
$u_{\varpi_{\ell_1}} \otimes \cdots \otimes u_{\varpi_{\ell_k}}$, where $u_{\varpi_i}$ is the highest 
weight element in $B(\varpi_i)$. As shown in~\cite{KN:1994},  the elements of $B(\om)$ can be labeled 
by tableaux of shape $\om$ in the alphabet $\{1,2,\ldots, n, \overline{n}, \ldots, \overline{1}\}$
for types $D_n$ and $C_n$ and the alphabet $\{1,2,\ldots, n, 0, \overline{n}, \ldots, \overline{1}\}$
for type $B_n$. For the explicit rules of type $D_n$, $B_n$, and $C_n$ tableaux we refer the reader 
to~\cite{KN:1994}; see also~\cite{HongKang:2002}.

\subsection{Definition of $\Bt^{r,s}$}
\label{ss:definition Btilde}

Let $\geh$ be of type $D_n^{(1)}$, $B_n^{(1)}$, or $A_{2n-1}^{(2)}$ with the 
underlying finite Lie algebra $\geh_0$ of type $X_n=D_n,B_n$, or $C_n$, respectively.
The combinatorial crystal $\Bt^{r,s}$ is defined as follows. As an $X_n$-crystal, $\Bt^{r,s}$
decomposes into the following irreducible components
\begin{equation} \label{eq:classical decomp}
	\Bt^{r,s} \cong \bigoplus_{\om} B(\om),
\end{equation}
for $1\le r\le n$ not a spin node. Here $B(\om)$ is the $X_n$-crystal of highest weight $\om$
and the sum runs over all dominant weights $\om$ that can be obtained from
$s\varpi_r$ by the removal of vertical dominoes, where $\varpi_i$ are the fundamental weights of 
$X_n$ as defined in section~\ref{ss:classical}.  The additional operators $\et{0}$ and $\ft{0}$
are defined as
\begin{equation} \label{eq:e0}
\begin{split}
\ft{0} &= \sigma \circ \ft{1} \circ \sigma,\\
\et{0} &= \sigma \circ \et{1} \circ \sigma,
\end{split}
\end{equation}
where $\sigma$ is the crystal analogue of the automorphism of the Dynkin diagram
that interchanges nodes 0 and 1. The involution $\sigma$ is defined in Definition~\ref{def:sigma}.

\subsection{Definition of $\sigma$}
\label{ss:sigma}

To define $\sigma$ we first need the notion of $\pm$ diagrams. A $\pm$ diagram $P$ of
shape $\La/\la$ is a sequence of partitions $\la\subset \mu \subset \La$ 
such that $\La/\mu$ and $\mu/\la$ are horizontal strips. We
depict this $\pm$ diagram by the skew tableau of shape $\La/\la$ in
which the cells of $\mu/\la$ are filled with the symbol $+$ and
those of $\La/\mu$ are filled with the symbol $-$. Write
$\La=\os(P)$ and $\la=\is(P)$ for the outer and inner shapes of the
$\pm$ diagram $P$. For type $A_{2n-1}^{(2)}$ and $r=n$, the inner shape $\lambda$
is not allowed to be of height $n$. When drawing partitions or tableaux, we use the French
convention where the parts are drawn in increasing order from top to bottom.

There is a bijection $\bij:P\mapsto b$ from $\pm$ diagrams $P$ of shape 
$\La/\la$ to the set of $X_{n-1}$-highest weight vectors $b$ of
$X_{n-1}$-weight $\la$ in $B_{X_n}(\La)$. Here $X_{n-1}$ is the subalgebra whose Dynkin 
diagram is obtained from that of $X_n$ by removing node $1$. There is a natural projection
of the weight lattices $\pi \colon P(X_n) \to P(X_{n-1})$, where $\pi(\alpha_i^{X_n})=
\alpha_{i-1}^{X_{n-1}}$ and $\pi(\varpi_i^{X_n})=\varpi_{i-1}^{X_{n-1}}$, and the partition
$\la$ is identified with the $X_{n-1}$ weights under $\pi$. We identify the Kashiwara operators
$\ft{i}^{X_{n-1}}$ with $\ft{i}^{X_n}$ under the embedding.

Explicitly the bijection $\bij$ is constructed as follows.
Define a string of operators 
$\ft{\aar} := \ft{a_1} \ft{a_2} \cdots \ft{a_\ell}$ such that $\bij(P) = \ft{\aar} u$, where $u$ is 
the highest weight vector in $B_{X_n}(\Lambda)$, where $\ft{i}$ is the Kashiwara
crystal operator corresponding to $f_i$.
Start with $\aar=()$. Scan the columns of $P$ from right
to left. For each column of $P$ for which a $+$ can be added, append
$(1,2, \ldots, h)$ to $\aar$, where $h$ is the height of the added $+$.
Next scan $P$ from left to right and for each column that contains a $-$ in $P$,
append to $\aar$ the string $(1,2,\ldots,n,n-2,n-3,\ldots, h)$ for type $D_n$,
$(1,2,\ldots,n-1,n,n,n-1,\ldots,h)$ for type $B_n$, and $(1,2,\ldots,n-1,n,n-1,\ldots,h)$ for
type $C_n$, where $h$ is the height of the $-$ in $P$. Note that for type $C_n$ the strings
$(1,2,\ldots,h)$ and $(1,2,\ldots,n-1,n,n-1,\ldots,h)$ are the same for $h=n$, which is why
empty columns of height $n$ are excluded for $\pm$ diagrams of type $A_{2n-1}^{(2)}$. 

By construction the automorphism $\sigma$ commutes with $\ft{i}$ and $\et{i}$ for $i=2,3,\ldots,n$. 
Hence it suffices to define $\sigma$ on $X_{n-1}$ highest weight elements. Because of the
bijection $\bij$ between $\pm$ diagrams and $X_{n-1}$-highest weight
elements, it suffices to define the map on $\pm$ diagrams.

Let $P$ be a $\pm$ diagram of shape $\La/\la$. Let $c_i=c_i(\la)$ be
the number of columns of height $i$ in $\la$ for all $1\le i<r$ with
$c_0=s-\la_1$. If $i\equiv r-1 \pmod{2}$, then in $P$, above each
column of $\la$ of height $i$, there must be a $+$ or a $-$.
Interchange the number of such $+$ and $-$ symbols. If $i\equiv r
\pmod{2}$, then in $P$, above each column of $\la$ of height $i$,
either there are no signs or a $\mp$ pair. Suppose there are $p_i$
$\mp$ pairs above the columns of height $i$. Change this to
$(c_i-p_i)$ $\mp$ pairs. The result is $\sigD(P)$, which has the
same inner shape $\la$ as $P$ but a possibly different outer shape.

\begin{definition} \label{def:sigma}
Let $b\in \Bt^{r,s}$ and $\et{\aar} := \et{a_1} \et{a_2} \cdots \et{a_\ell}$ be such that
$\et{\aar}(b)$ is a $X_{n-1}$ highest weight crystal element. Define 
$\ft{\aal}:= \ft{a_\ell} \ft{a_{\ell-1}} \cdots \ft{a_1}$. Then 
\begin{equation} \label{eq:def sigma}
\sigma(b) := \ft{\aal} \circ \bij \circ \sigD \circ \bij^{-1} \circ \et{\aar}(b).
\end{equation}
\end{definition}
It was shown in~\cite{S:2007} that $\Bt^{r,s}$ is regular.

\subsection{Properties of $\Bt^{r,s}$}

For the proof of uniqueness we will require the action of $\et{1}$ on $X_{n-2}$ highest
weight elements, where $X_{n-2}$ is the Dynkin diagram obtained by removing nodes
1 and 2 from $X_n$. As we have seen in section~\ref{ss:sigma}, the $X_{n-1}$-highest weight 
elements in the branching $X_n\to X_{n-1}$ can be described by $\pm$ diagrams. Similarly the 
$X_{n-2}$-highest weight elements in the branching $X_{n-1}\to X_{n-2}$ can be described 
by $\pm$ diagrams. Hence each $X_{n-2}$-highest weight vector is uniquely determined by a 
pair of $\pm$ diagrams $(P,p)$ such that $\is(P)=\os(p)$. The diagram $P$ specifies the 
$X_{n-1}$-component $B_{X_{n-1}}(\is(P))$ in $B_{X_n}(\os(P))$, and $p$ specifies the 
$X_{n-2}$ component inside $B_{X_{n-1}}(\is(P))$. Let $\Upsilon$ denote the map $(P,p) \mapsto b$ 
from a pair of $\pm$ diagrams to a $X_{n-2}$ highest weight vector.

To describe the action of $\et{1}$ on an $X_{n-2}$ highest weight element or by $\Upsilon$ 
equivalently on $(P,p)$ perform the following algorithm:
\begin{enumerate}
\item Successively run through all $+$ in $p$ from left to right and, if possible, pair it with 
the leftmost yet unpaired $+$ in $P$ weakly to the left of it.
\item Successively run through all $-$ in $p$ from left to right and, if possible, pair it with
the rightmost  yet unpaired $-$ in $P$ weakly to the left.
\item Successively run through all yet unpaired $+$ in $p$ from left to right and, if possible,
pair it with the leftmost yet unpaired $-$ in $p$.
\end{enumerate}

\begin{lemma} \cite[Lemma 5.1]{S:2007} \label{lem:e1}
If there is an unpaired $+$ in $p$,  $\et{1}$ moves the rightmost unpaired $+$ in $p$ to $P$. 
Otherwise, if there is an unpaired $-$ in $P$, $\et{1}$ moves the leftmost unpaired $-$ in $P$ to $p$.
Otherwise $\et{1}$ annihilates $(P,p)$.
\end{lemma}

In this paper, we will only require the case of Lemma~\ref{lem:e1} when a $-$ from $P$ moves
to $p$. Schematically, if a $-$ from a $\mp$ pair in $P$ moves to $p$, then the following happens
\begin{equation*}
\young(\mb\bluem-,\mb++) \; \raisebox{.3cm}{$\mapsto$} \; \young(\mb+-,\mb\redm+) 
\quad \raisebox{0.3cm}{or} \quad
\young(--\bluem-,\mb\mb++) \; \raisebox{.3cm}{$\mapsto$} \; \young(+---,\mb\mb\redm+)
\; \raisebox{.3cm}{,}
\end{equation*}
where the blue minus is the minus in $P$ that is being moved and the red minus is the new minus
in $p$. Similarly, schematically if a $-$ not part of a $\mp$ pair in $P$ moves to $p$, then
\begin{equation*}
\young(\mb\bluem--,\mb\mb\mb\mb) \; \raisebox{.3cm}{$\mapsto$} \; \young(\mb\redm--,\mb\mb\mb\mb)
\quad \raisebox{0.3cm}{or} \quad
\young(++\bluem--,\mb\mb\mb\mb\mb) \; \raisebox{.3cm}{$\mapsto$} \; 
\young(\redm++--,\mb\mb\mb\mb\mb) \; \raisebox{.3cm}{.}
\end{equation*}

For any $b\in \Bt^{r,s}$, let $\is(b)$ be the inner shape of the $\pm$ diagram corresponding
to the $X_{n-1}$ highest weight element in the component of $b$. Furthermore recall that
$\Bt^{r,s}$ is regular, so that in particular $\et{0}$ and $\et{1}$ commute. 
We can now state the lemma needed in the next section.

\begin{lemma} \label{lem:containment}
Let $b\in\tilde{B}^{r,s}$ be an $X_{n-2}$ highest weight vector corresponding under 
$\Upsilon$ to the tuple of $\pm$ diagrams $(P,p)$ where $\is(p)=\os(p)$. 
Assume that $\ve_0(b),\ve_1(b)>0$. Then $\is(b)$ is strictly contained in 
$\is(\et{0}(b))$, $\is(\et{1}(b))$, and $\is(\et{0}\et{1}(b))$.
\end{lemma}

\begin{proof}
By assumption $p$ does not contain any $-$ and $\et{1}$ is defined. Hence $\et{1}$
moves a $-$ in $P$ to $p$. This implies that the inner shape of $b$ is strictly
contained in the inner shape of $\et{1}(b)$.

The involution $\sigma$ does not change the inner shape of $b$ (only the outer shape).
By the same arguments as before, the inner shape of $b$ is strictly contained in
the inner shape of $\et{1} \sigma(b)$. Since $\sigma$ does not change the inner shape, this 
is still true for $\et{0}(b) = \sigma \et{1} \sigma(b)$.

Now let us consider $\et{0} \et{1}(b)$. For the change in inner shape
we only need to consider $\et{1} \sigma \et{1}(b)$, since the last $\sigma$
does not change the inner shape. By the same arguments as before,  $\et{1}$
moves a $-$ from $P$ to $p$ and $\sigma$ does not change the inner shape. The next
$\et{1}$ will move another $-$ in $\sigma \et{1}(b)$ to $p$. Hence $p$ will
have grown by two $-$, so that the inner shape of $\et{1} \sigma \et{1}(b)$
is increased by two boxes.
\end{proof}

\section{Equivalence of $B^{r,s}$ and $\Bt^{r,s}$ of type $D_n^{(1)}$, $B_n^{(1)}$,
and $A_{2n-1}^{(2)}$}
\label{sec:equivalence}

In this section all crystals are of type $D_n^{(1)}$, $B_n^{(1)}$, or $A_{2n-1}^{(2)}$ with
corresponding classical subalgebra of type $X_n=D_n,B_n,C_n$, respectively.

Let $B$ and $B'$ be regular crystals of type $D_n^{(1)}$, $B_n^{(1)}$, or $A_{2n-1}^{(2)}$ with index 
set $I=\{0,1,2,\ldots,n\}$. We say that $B\simeq B'$ is an isomorphism of
$J$-crystals if $B$ and $B'$ agree as sets and all arrows colored $i\in J$ are the same. 

\begin{proposition} \label{prop:Psi}
Suppose that there exist two isomorphisms
\begin{equation*}
\begin{split}
	\Psi_0 : \tilde{B}^{r,s} \simeq B & \qquad \text{as an isomorphism of $\{1,2,\ldots,n\}$-crystals}\\
	\Psi_1 : \tilde{B}^{r,s} \simeq B & \qquad \text{as an isomorphism of $\{0,2,\ldots,n\}$-crystals.}
\end{split}
\end{equation*}
Then $\Psi_0(b) = \Psi_1(b)$ for all $b\in \tilde{B}^{r,s}$ and hence there exists an $I$-crystal
isomorphism $\Psi:\tilde{B}^{r,s} \simeq B$.
\end{proposition}

\begin{remark}
Note that $\Psi_0$ and $\Psi_1$ preserve weights, that is,
$\wt(b)=\wt(\Psi_0(b))=\wt(\Psi_1(b))$ for all $b\in \tilde{B}^{r,s}$.
This is due to the fact that if all but one coefficient $m_j$ are known for a weight
$\Lambda=\sum_{j=0}^n m_j \La_j$, then the missing $m_j$ is also determined by the
level 0 condition.
\end{remark}

\begin{proof}
If $\Psi_0(b)=\Psi_1(b)$ for a $b$ in a given $X_{n-1}$-component $C$, then $\Psi_0(b')=\Psi_1(b')$ 
for all $b'\in C$ since $\et{i}\Psi_0(b')=\Psi_0(\et{i}b')$ and $\et{i}\Psi_1(b')=\Psi_1(\et{i}b')$ for 
$i\in J=\{2,3,\ldots,n\}$. Hence it suffices to prove $\Psi_0(b)=\Psi_1(b)$ for only one element
$b$ in each $X_{n-1}$-component $C$. We are going to establish the theorem for $b$ corresponding to
the pairs of $\pm$ diagrams $(P,p)$ where $\is(p)=\os(p)$. Note that this is an $X_{n-2}$-highest
weight vector, but not necessarily an $X_{n-1}$-highest weight vector.

We proceed by induction on $\is(b)$ by containment. First suppose that both $\ve_0(b),\ve_1(b)>0$.
By Lemma~\ref{lem:containment}, the inner shape of $\et{0} \et{1}b$, $\et{0}b$, and $\et{1}b$ is bigger
than the inner shape of $b$, so that by induction hypothesis
$\Psi_0(\et{0} \et{1}b)=\Psi_1(\et{0} \et{1}b)$, $\Psi_0(\et{0}b)=\Psi_1(\et{0}b)$, and 
$\Psi_0(\et{1}b)=\Psi_1(\et{1}b)$. Therefore we obtain
\begin{multline*}
	\et{0} \et{1} \Psi_0(b) = \et{0} \Psi_0(\et{1} b) = \et{0} \Psi_1(\et{1} b) = \Psi_1(\et{0} \et{1} b)
	= \Psi_0(\et{0} \et{1} b)\\
	 = \et{1} \Psi_0(\et{0} b) = \et{1} \Psi_1(\et{0} b) = \et{1} \et{0} \Psi_1(b).
\end{multline*}
This implies that $\Psi_0(b)=\Psi_1(b)$.

Next we need to consider the cases when $\ve_0(b)=0$ or $\ve_1(b)=0$, which comprises the base
case of the induction. Let us first treat the case $\ve_1(b)=0$.
Recall that $\is(p)=\os(p)$ so that $p$ contains only empty columns.
Hence it follows from the description of the action of $\et{1}$ of Lemma~\ref{lem:e1},
that $\ve_1(b)=0$ if and only if $P$ consists only of empty columns or columns containing $+$. 
\begin{quote}
\textbf{Claim.} $\Psi_0(b)=\Psi_1(b)$ for all $b$ corresponding to the pair of $\pm$ diagrams
$(P,p)$ where $P$ contains only empty columns and columns with $+$, and $\is(p)=\os(p)$.
\end{quote}
The claim is proved by induction on $k$, which is defined to be the number of empty columns
in $P$ of height strictly smaller than $r$. For $k=0$ the claim is true by weight considerations.
Now assume the claim is true for all $0\le k'<k$ and we will establish the claim for $k$.
Suppose that $\Psi_1(b)=\Psi_0(\tilde{b})$ where $\tilde{b}\neq b$. By weight considerations 
$\tilde{b}$ must correspond to a pair of $\pm$ diagrams $(\tilde{P},p)$, where $\tilde{P}$
has the same columns containing $+$ as $P$, but some of the empty columns of $P$ of
height $h$ strictly smaller than $r$ could be replaced by columns of height $h+2$ containing
$\mp$. Denote by $k_+$ the number of columns of $P$ containing $+$. Then
\begin{equation*}
	m := \ve_0(b) = k_+ + k,
\end{equation*}
since under $\sigma$ all empty columns in $P$ become columns with $\pm$ and columns
containing $+$ become columns with $-$. By Lemma~\ref{lem:e1}, then $\et{1}$ acts on
$(\sigD(P),p)$ as often as there are minus signs in $\sigD(P)$, which is $k_++k$.
Set $\hat{b}=\et{1}^a\tilde{b}$, where $a>0$ is the number of columns in $\tilde{P}$ containing
$\mp$. If $(\hat{P},\hat{p})$ denotes the tuple of $\pm$ diagrams associated
to $\hat{b}$, then compared to $(\tilde{P},p)$ all $-$ from the $\mp$ pairs in $\tilde{P}$ moved to $p$.
Note that $\hat{P}$ has only $k-a<k$ empty columns of height less than $r$, so that by induction
hypothesis $\Psi_0(\hat{b})=\Psi_1(\hat{b})$. Hence
\begin{equation} \label{eq:eq string}
	\Psi_1(b) = \Psi_0(\tilde{b}) = \Psi_0(\ft{1}^a \hat{b}) = \ft{1}^a \Psi_0(\hat{b}) 
	= \ft{1}^a \Psi_1(\hat{b}).
\end{equation}
Note that 
\begin{equation*}
	\ve_0(\hat{b}) = \ve_0(\tilde{b}) = m-a<m.
\end{equation*}
Hence
\begin{equation*}
\begin{split}
	&\et{0}^m \Psi_1(b) = \Psi_1(\et{0}^m b) \neq 0\\
	\text{but} \qquad 
	&\et{0}^m \ft{1}^a \Psi_1(\hat{b}) = \ft{1}^a \Psi_1(\et{0}^m \hat{b}) = 0
\end{split}
\end{equation*}
which contradicts~\eqref{eq:eq string}. This implies that we must have $\tilde{b}=b$ proving the claim.

The case $\varepsilon_0(b)=0$ can be proven in a similar fashion to the case $\varepsilon_1(b)=0$. 
Using the explicit action of $\sigD$ on $P$ and Lemma~\ref{lem:e1}, it follows that $\varepsilon_0(b)=0$
if and only if $P$ consists only of columns containing $-$ or $\mp$ pairs. 
\begin{quote}
\textbf{Claim.} $\Psi_0(b)=\Psi_1(b)$ for all $b$ corresponding to the pair of $\pm$ diagrams
$(P,p)$ where $P$ contains only columns with $-$ and columns with $\mp$ pairs, and $\is(p)=\os(p)$.
\end{quote}
By induction on the number of $\mp$ pairs in $P$, this claim can be proven similarly as before
(using the fact that $\sigD$ changes columns with $-$ into columns with $+$ and columns
with $\mp$ pairs into empty columns).
\end{proof}

\begin{proof}[Proof of Theorem~\ref{thm:equiv}]
Both crystals $B^{r,s}$ and $\Bt^{r,s}$ have the same classical 
decomposition~\eqref{eq:classical decomp} as $X_n$ crystals with index set $\{1,2,\ldots,n\}$ 
and $\{0,2,3,\ldots,n\}$ by Proposition~\ref{prop:classical decomp}. Hence there exist
crystal isomorphisms $\Psi_0$ and $\Psi_1$. By Proposition~\ref{prop:Psi} there exists
an $I$-isomorphism $\Psi: \Bt^{r,s} \cong B^{r,s}$ which proves the theorem.
\end{proof}

\appendix

\section{Erratum} \label{erratum}
Here we would like to correct some errors and omissions in our paper, that we noticed after
publication.
\begin{enumerate}
\item In Table 1, node $n$ for type $B_n^{(1)}$ should not be filled.
Also, the terminology "spin node" as used in Section 4.2 is misleading.
In~\cite{FOS:2008} we use the terminology "exceptional node" instead.
\item In Section 4.2, the decomposition of $W_s^{(n)}$ for $B_n^{(1)}$
as a $U_q(\mathfrak{g}_0)$-module should be given by Eq. (4.1), where we identify
$\varpi_n$ with a column of height $n$ and of width $1/2$, and $\omega$ runs over all 
partitions that can be obtained from the $n\times(s/2)$ rectangle by removing vertical
dominoes.
\item The proof of Proposition 6.1 does not apply to the case of type 
$A_{2n-1}^{(2)}$ and columns of height $n$ (since the inner $\pm$-diagram $p$ is not
allowed to have empty columns). See the proof of~\cite[Theorem 5.1]{FOS:2008}
for this case.
 \end{enumerate}

The first paragraph of Section 4.3 needs to be extended to the case $r=n$ for 
$B_n^{(1)}$, which is done below.

\subsection{Calculation of prepolarization: $B_n^{(1)},r=n$ case}
Let $n'=[n/2]$. 
Let $\cb=(c_1,c_2,\ldots,c_{n'})$ be a sequence of integers such that 
$s/2\ge c_1\ge c_2\ge\cd\ge c_{n'}\ge0$. For such $\cb$ we define a vector $u_m$ ($0\le m\le n'$) 
in $W_s^{(n)}$ inductively by
\[
u_m=(e_{n-2m}^{(c_m)}\cd e_2^{(c_m)}e_1^{(c_m)})(e_{n-2m+1}^{(c_m)}\cd e_3^{(c_m)}e_2^{(c_m)})
e_0^{(c_m)}u_{m-1},
\]
where $u_0$ is the vector in (iii) of Proposition \ref{prop:KRassump}. 
Set $u(\cb)=u_{n'}$. The weight of $u(\cb)$ is given by
\[
\la(\cb)=\sum_{j=0}^{n'}(c_j-c_{j+1})(1+\delta_{j0})\varpi_{n-2j},
\]
where we have set $c_0=s/2,c_{n'+1}=0$, and $\varpi_0$ should be understood as $0$. 
$\la(\cb)$ represents all $\omega$ in~\eqref{branching} 
when $\cb$ runs over all possible sequences.
The following proposition calculates values of the prepolarization $(\;,\;)$ on $W_s^{(n)}$.

\newpage

\begin{proposition}\mbox{}
\begin{itemize}
\item[(1)] ${\displaystyle (u(\cb),u(\cb))=\prod_{j=1}^{n'}q^{c_j(s-c_j)}{s\brack c_j}}$,
\item[(2)] $(e_ju(\cb),e_ju(\cb))=0$ unless $n-j\in 2\Z_{\ge0}$. If $n-j\in 2\Z_{\ge0}$, then 
	setting $p=(n-j)/2+1$, $(e_ju(\cb),e_ju(\cb))$ is given by
\[
\prod_{j=1}^{n'}q^{(c_j-\delta_{j,p})(s-c_j)}
{s-\delta_{j,p}\brack c_j-\delta_{j,p}}
\times\left\{
\begin{array}{ll}
q_n^{s-1}[s]_n
\quad&\text{if }p=1,\\
q^{s-c_{p-1}-1}[s-c_{p-1}]
\quad&\text{if }p>1.
\end{array}
\right.
\]
\end{itemize}
\end{proposition}

\end{document}